% this is, as of today (January 27, 2007), the final
% version of our paper entitled
% Weighted Poincar\'e inequality and rigidity of complete manifolds

\input amstex
\documentstyle{amsppt}
%\magnification=1200
\vsize=7.4in
\hsize=5.5in
\nologo

%\noindent
%$\underline{ Ann. Sc. Ec. Norm. Sup.}$  $4^e$ serie, t. 0, 2007, p. 1--62
%\bigskip

\def\d{\Delta}
\def\p{\partial}
\def\n{\nabla}
\def\la{\langle}
\def\ra{\rangle}
\def\br{\Bbb R}     
\def\-{\setminus}

\topmatter
\title
Weighted Poincar\'e inequality and rigidity of complete manifolds
\endtitle

\rightheadtext{Weighted Poincar\'e inequality and rigidity}
\leftheadtext{Peter Li and Jiaping Wang}
\author
Peter Li and Jiaping Wang
\endauthor
\affil
University of California, Irvine\\
University of Minnesota
\endaffil

\footnotetext" "{The first author was partially supported by NSF grant DMS-0503735. 
The second author was partially supported by NSF grant DMS-0404817.}

\address
Department of Mathematics, University of California,
Irvine, CA 92697-3875
\endaddress

\email
pli\@math.uci.edu
\endemail

\address
School of Mathematics, University of Minnesota,
Minneapolis, MN 55455
\endaddress

\email
jiaping\@math.umn.edu
\endemail

\date
%Jan 27, 2007
\enddate

\abstract
We prove structure theorems for complete manifolds satisfying both the Ricci curvature lower bound
and the weighted Poincar\'e inequality. In the process, a sharp decay estimate for the minimal positive
Green's function is obtained. This estimate only depends on the weight function of the Poincar\'e inequality,
and yields a criterion of parabolicity of connected components at infinity in terms of the weight function. 

\medskip

\centerline {\bf In\'egalit\'es de Poincar\'e \`a poids et rigidit\'e des vari\'et\'es compl\`etes}

\smallskip

\noindent{R\'esum\'e:
Nous prouvons des th\'eor\`emes de structure pour des vari\'et\'es compl\`etes telles que la courbure de Ricci 
soit minor\'ee, 
et satisfaisant l'in\'egalit\'e de Poincar\'e \`a poids.
Nous obtenons une  estimation optimale de la d\'ecroissance de la fonction de Green positive et minimale. 
Cette estimation, qui d\'epend seulement du poid de la fonction
dans l'in\'egalit\'e de Poincar\'e, produit un crit\'ere de parabolicit\'e de composantes reliées à l'infini 
utilisant le poid de la fonction.}
\endabstract

\endtopmatter

\document

\subheading
{Table of Contents}

\noindent
\S0 Introduction \newline
\noindent
\S1 Weighted Poincar\'e Inequality \newline
\noindent
\S2 Decay Estimate \newline
\noindent
\S3 Geometric Conditions for Parabolicity and Nonparabolicity \newline
\noindent
\S4 Improved Bochner Formula and Metric Rigidity \newline
\noindent
\S5 Rigidity and Nonparabolic Ends \newline
\noindent
\S6 Warped Product Metrics \newline
\noindent
\S7 Parabolic Ends \newline
\noindent
\S8 Nonexistance Results for Parabolic Ends

\heading
\S0 Introduction
\endheading

Understanding the relations among the curvature, the topology and the function theory is a 
central theme in Riemannian geometry.
Typically, one assumes the curvature to be bounded by a constant so that the comparison theorems become available. 
The main focus of 
the current paper, however, is to go beyond this realm to consider manifolds with Ricci curvature bounded below 
by a function.
We will establish some structure theorems for such manifolds satisfying the following Poincar\'e type inequality.   

\proclaim{Definition 0.1} Let $M^n$ be an $n$-dimensional complete Riemannian manifold.  
We say that $M$ satisfies a {\rm weighted Poincar\'e inequality} 
with a nonnegative weight function $\rho(x)$, if the inequality
$$
\int_M \rho(x)\, \phi^2(x)\,dV \le \int_M |\n \phi |^2\, dV
$$
is valid for all compactly supported smooth function $\phi \in C_c^{\infty}(M)$.
\endproclaim  

\proclaim{Definition 0.2} Let $M^n$ be an $n$-dimensional complete Riemannian manifold.  We say that $M$ 
has property ($\Cal P_\rho$) if a 
weighted Poincar\'e inequality is valid on $M$ with some nonnegative weight function $\rho$.  Moreover, 
the $\rho$-metric, defined by
$$
ds^2_\rho = \rho\, ds^2_M
$$
is complete.
\endproclaim

Let $\lambda_1(M)$ denote the greatest lower bound of the spectrum of the Laplacian acting on $L^2$ functions. Then  
the variational principle for $\lambda_1(M)$ asserts the validity of the Poincar\'e inequality, i.e., 
$$
\lambda_1(M) \int_M \phi^2 \le \int_M |\n \phi|^2
$$
for all compactly supported functions $\phi \in C_c^{\infty}(M).$ 
Obviously, $M$ has property ($\Cal P_\rho$) with weight function $\rho = \lambda_1(M)$ if $\lambda_1(M)>0.$  
Hence the notion of property ($\Cal P_\rho$) can be viewed as a generalization of the assumption $\lambda_1(M) >0.$

We would like to point out that the idea of considering $ds^2_\rho$ was first used by 
Agmon \cite{A} in his study of eigenfunctions for the Schr\"odinger operators. Indeed, we will employ
some of the arguments from \cite{A} in this paper. We also remark that the
weighted Poincar\'e inequalities in various forms have appeared in many important issues of analysis 
and mathematical physics. 
In the interesting papers \cite{F-P1} and \cite{F-P2}, Fefferman and Phong have considered the more general
weighted Sobolev type inequalities for pseudodifferential operators. 
   
In \S1, we will demonstrate that a complete manifold is nonparabolic if and only if it satisfies 
a weighted Poincar\'e inequality 
with some weight function $\rho.$  Moreover, many nonparabolic manifolds satisfy
property ($\Cal P_\rho$) and we will  provide some 
systematic ways to find a weight function.  However, we would like to point out that the weight 
function is obviously not unique. 
  
It turned out that some of the crucial estimates we developed in \cite{L-W1} can be generalized to give a sharper version of 
Agmon's estimate for the Schr\"odinger operators developed in \cite{A}.  We believe that these estimates are interesting in their own 
rights and the sharp 
form will find more geometric applications in our further investigation.    
The proof of this decay estimate will be given in \S2.  

In  \S3, we will give geometric conditions (involving $\rho$) for an end being nonparabolic or parabolic using the decay estimate obtained in \S2.  
The  conditions are parallel to what we have established for the case when $\rho= \lambda_1(M)$ in \cite{L-W1}.  In \S4, we recall a generalized Bochner 
formula for the gradient of a harmonic function.  The equality case for this inequality will be discussed.  Note that the generalized Bochner 
formula was first used by Yau in \cite{Y1} and the equality case was also used previously 
in \cite{W} and \cite{L-W1}.  

In \S5, we give the proof of a structure theorem (Theorem 5.2) for manifolds with property ($\Cal P_\rho$).
The Ricci curvature is assumed to satisfy the inequality
$$
\text{Ric}_M(x) \ge -\frac{n-1}{n-2}\, \rho(x)
$$
for all $x \in M.$
In this theorem, a growth assumption is needed for the weight function $\rho$ in terms of the $\rho$-distance.  This allows to account for those 
functions $\rho$ that has different growth rate with respect to the background distance in different directions.  In any case, if $\rho$ is 
bounded then the growth assumption is fulfilled.  When the dimension of $M$ is at least 4, the assumption on $\rho$ is rather mild.  
In particular, if $\rho$ is a non-decreasing function of the distance to a compact set (with respect to the background metric) then 
the growth assumption is automatically satisfied.  

In \S6, we will study when a warped-product situation will actually occur for specific choices of $\rho.$  When $\rho$ is an increasing function, 
we give examples of the warped-product scenario.  On the other hand, we will also prove in Theorem 6.3 that when $\liminf_{r \to \infty} \rho(r) = 0$, 
the warped-product scenario does not exist. Let us summarize the results in \S5 and \S6 in the following Theorem. We denote by $S(R)$ the maximum value of $\sqrt {\rho}$
over the geodesic ball of radius $R$ with respect to the $\rho$-metric centered at a fixed point $p.$ 

\proclaim{Theorem A} Let $M^n$ be a complete manifold with dimension $n \ge 4.$  Assume that $M$ satisfies property ($\Cal P_\rho$) for some 
nonzero weight function $\rho \ge 0.$  Suppose the Ricci curvature of $M$ satisfies the lower bound
$$
\text{Ric}_M(x) \ge - \frac{n-1}{n-2}\, \rho(x)
$$
for all $x \in M.$  If $\rho$ satisfies the growth estimate

$$
\liminf_{R\to \infty} S(R)\, \exp \left(-\frac{n-3}{n-2}\, R \right) =0,
$$
then either
\roster
\item $M$ has only one nonparabolic end; or
\item $M$ has two nonparabolic ends and is given by $M = \Bbb R \times N$ with the warped product metric
$$
ds_M^2 = dt^2 + \eta^2(t)\,ds_N^2,
$$
for some positive function $\eta(t)$, and some compact manifold $N.$ Moreover, $\rho(t)$ is a function of $t$ alone satisfying
$$
(n-2)\eta''\, \eta^{-1} = \rho
$$
and
$$
\liminf_{x \to \infty}\rho(x) >0.
$$
\endroster
\endproclaim

We will state the case $n=3$ separately since the curvature condition for counting parabolic and nonparabolic ends are the same in this case.  
Again, the condition on the growth rate of $\rho$ is given with respect to the $\rho$-distance.  When restricted to those weight functions
that are non-decreasing functions of the distance to the compact set the growth assumption is simply subexponential growth with respect to 
the background distance function.  Combining with the results in \S7 for parabolic ends, the 3-dimension case can be stated as follows:

\proclaim{Theorem B} Let $M^3$ be a complete manifold of dimension  3.  Assume that $M$ satisfies property ($\Cal P_\rho$) for 
some nonzero weight function $\rho \ge 0.$  Suppose the Ricci curvature of $M$ satisfies the lower bound
$$
\text{Ric}_M(x) \ge - 2\, \rho(x)
$$
for all $x \in M.$  If $\rho$ satisfies the growth estimate
$$
\liminf_{R\to \infty} S(R)\, R^{-1}=0,
$$
then either
\roster
\item $M$ has only one end; 
\item $M$ has two nonparabolic ends and is given by $M = \Bbb R \times N$ with the warped product metric
$$
ds_M^2 = dt^2 + \eta^2(t)\,ds_N^2,
$$
for some positive function $\eta(t)$, and some compact manifold $N.$ Moreover, $\rho(t)$ is a function of $t$ alone satisfying
$$
\eta''\, \eta^{-1} = \rho
$$
and
$$
\liminf_{x \to \infty} \rho(x) >0; \text{or}
$$
\item $M$ has one parabolic end and one nonparabolic end and is given by $M = \Bbb R \times N$ with the warped product metric
$$
ds_M^2 = dt^2 + \eta^2(t)\,ds_N^2,
$$
for some positive function $\eta(t)$, and some compact manifold $N.$ Moreover, $\rho(t)$ is a function of $t$ alone satisfying
$$
\eta''\, \eta^{-1} = \rho
$$
and
$$
\liminf_{x \to \infty} \rho(x) >0 \  \text{on the nonparabolic end.}
$$
\endroster
\endproclaim

Also, for dimension $n \ge 4$ we proved the following:

\proclaim{Theorem C} Let $M^n$ be a complete manifold of dimension $n \ge 4$ with property ($\Cal P_\rho$).  Suppose the Ricci curvature of $M$ satisfies the lower bound
$$
\text{Ric}_M (x) \ge -\frac{4}{n-1}\, \rho(x)
$$
for all $x \in M.$  If $\rho$ satisfies the property that
$$
\lim_{x\to \infty} \rho(x)=0,
$$
then $M$ has only one end.
\endproclaim

Finally, in the last section, \S8, we prove a nonexistence result indicating that for a large class of weight functions $\rho$, 
namely when $\rho$ is a function of the distance and satisfying $(\rho^{-\frac 14})''(r)\ge 0$ for $r$ sufficiently large, there  
does not exist a manifold with property ($\Cal P_\rho$) satisfying
$$
\text{Ric}_M (x) \ge - \frac 4{n-1}\, \rho(x).
$$
We also proved a theorem restricting the behavior of the warped product.

As we have pointed out earlier, if $\lambda_1(M)>0,$ then one may take $\rho=\lambda_1(M).$ This special case of
Theorem A and Theorem B has been the subject of our earlier work \cite{L-W1} and \cite{L-W2}. The results generalized 
the work of Witten-Yau \cite{W-Y}, Cai-Galloway \cite{C-G}, and Wang \cite{W}, on conformally compact manifolds. Also,
a result in the similar spirit of Theorem C is available for this special case (see \cite{L-W2}). However, it remains
open to deal with more general functions $\rho$ in Theorem C.

We would also like to point out that a similar theory was proposed in \cite{L-W3} for K\"ahler manifolds where the assumption is on 
the holomorphic bisectional curvature instead.

The authors would like to thank S.T. Yau for his interest in this work and his encouragement for  us to continue our effort in considering 
the general case when $\rho$ is not only a function of the distance.

\heading
\S1 Weighted Poincar\'e Inequality
\endheading

In this section, we will show that it is not difficult to find a weight function for most manifolds.
The following proposition gives  a convenient way to construct a weight function.   The arugment uses  
a modified version of the Barta inequality for the first eigenvalue.  Using this method, we will show that a manifold is nonparabolic 
if and only if there exists a non-trivial weight function. We will also give examples of using other methods to obtain a weight function for various manifolds.

\proclaim{Proposition 1.1} Let $M$ be a complete Riemannian manifold.  If there 
exists a nonnegative function $h$ defined on $M$, that is not identically 0, satisfying
$$
\d h(x)\le -\rho(x)\, h(x),
$$
then the weighted Poincar\'e inequality
$$
\int_M \rho(x)\, \phi^2(x) \le \int_M |\n \phi|^2(x)
$$
must be valid for all compactly supported smooth function $\phi \in C_c^{\infty}(M).$
\endproclaim

\demo{Proof}Let $D\subset M$ be a smooth compact subdomain of $M.$ Let us denote 
$\lambda_1(\rho, D)$ to be the first Dirichlet eigenvalue on $D$  for the operator
$$
\d + \rho(x).
$$
Let $u$ be the first eigenfunction satisfying
$$
\d u(x)  + \rho(x)\, u(x)= - \lambda_1(\rho, D) \,u(x) \qquad \text{on} \qquad D
$$
and
$$
u(x) = 0 \qquad \text{on} \qquad \partial D.
$$
We may assume that $u \ge 0$ on $D,$ and the regularity of $u$ asserts that $u>0$ in 
the interior of $D.$  Integration by parts yields
$$
\split
\int_D u\,\Delta h - \int_D h\,\Delta u 
& = \int_{\p D} u \, \frac{\partial h}{\partial \nu} - \int_{\p D} h \, \frac{\partial 
u}{\partial \nu}\\
& \ge 0,
\endsplit\tag 1.1
$$
where $\nu$ is the outward unit normal of $\p D$. On the other hand, the assumption on $h$ implies that
$$
u\, \d h - h\, \d u \le   \lambda_1(\rho, D)\,u\, h.
$$
Since both $u>0$ and $h$ are not identically 0, this combining with \thetag {1.1} implies that $\lambda_1(\rho, D) \ge 0.$    
In particular, the variational characterization of $\lambda_1(\rho, D)$ implies that 
$$
\split
0 &\le 
\lambda_1(\rho, D) \int_D \phi^2(x) \\
&\le \int_D |\n \phi|^2 - \int_D \rho(x)\, \phi^2(x)
\endsplit
$$
for all $\phi$ with support in $D.$
Since $D$ is arbitrary, this implies the weighted Poincar\'e inequality.\qed
\enddemo

Let us now quickly recall the definition of parabolicity.   Full detail on the discussion of parabolicity can be found in \cite{L2}.

\proclaim{Definition 1.2} A complete manifold $M$ is said to be nonparabolic if there exists a symmetric positive Green's function $G(x, y)$ 
for the Laplacian acting on $L^2$ functions.  Otherwise, we say that $M$ is parabolic.
\endproclaim

Similarly, the notion of parabolicity is also valid when localized at an end of a manifold.  We recall (see \cite{L2})  
that an end is simply an unbounded component of $M \- \Omega$, where $\Omega$ is a compact smooth domain of $M$.  

\proclaim{Definition 1.3} Let $E$ be an end of a complete manifold $M$.  We say that $E$ is nonparabolic if there exists a symmetic, 
positive, Green's function $G(x, y)$ for the Laplacian acting on $L^2$ functions with Neumann boundary condition on $\p E.$  
Otherwise, we say that $E$ is a parabolic end.
\endproclaim

Note that (see \cite{L2}) a complete manifold is nonparabolic if and only if $M$ has at least one nonparabolic end.   More importantly, 
it is possible for a nonparabolic manifold to have many parabolic ends.  It is also known that \cite{L2} a manifold is parabolic if we 
consider the sequence of harmonic functions $\{f_i\}$ defined on $B(R_i)\- B(R_0)$ satisfying
$$
\d f_i = 0 \qquad \text{on} \qquad B(R_i)\- B(R_0)\tag 1.2
$$
with boundary conditions
$$
f_i = \left\{ \aligned 1& \qquad \text{on} \qquad \p B(R_0)\\
0& \qquad \text{on} \qquad \p B(R_i),\endaligned \right.\tag 1.3
$$
then they converge to the constant function $f = 1$ defined on $M \- B(R_0)$ as $R_i \to \infty$ for any fixed $R_0.$

\proclaim{Corollary 1.4} Let $M$ be a complete  nonparabolic manifold and $G(p, x)$ be the minimal positive Green's function defined on $M$ 
with a pole at the point $p \in M.$  Then $M$ satisfies the weighted Poincar\'e inequality with the weight function $\rho$ given by
$$
\rho (x) = \frac{|\n G(p, x)|^2}{4G^2(p, x)}.
$$ 
Conversely, if a non-zero weight function $\rho \ge 0$ exists, then $M$ must be nonparabolic.
\endproclaim

\demo{Proof} Let us first assume that $M$ is nonparabolic and hence a positive symmetric Green's function, $G(x, y),$ for the Laplacian exists.
Let $p \in M$ be a fixed point and
$$
g_a(x) =\min \{ a, G(p, x)\}.
$$
Then $g_a$ is a superharmonic function defined on $M.$  A direct computation yields that
$$
\d g_a^{\frac 12} \le - \frac{|\n g_a|^2}{4g_a^2}\, g_a^{\frac 12}
$$
on $M.$  Hence Proposition 1.1 asserts that the weighted Poincar\'e inequality is valid with
$$
\rho =\frac{|\n g_a|^2}{4g_a^2}.
$$
Letting $a \to \infty,$ we conclude that we can take 
$$
\rho (x) = \frac{|\n G(p, x)|^2}{4G^2(p, x)}.
$$
This proves the first part of the corollary.

Conversely, let us assume that the weighted Poincar\'e inequality is valid for a non-zero weight function $\rho \ge 0.$ 
Assuming the contrary that $M$ is parabolic, we will find a contradiction.  Indeed, if $M$ is parabolic then let us consider the 
sequence of compactly supported functions
$$
\phi_i = \left\{ \aligned 1& \qquad \text{on} \qquad B(R_0)\\
f_i& \qquad \text{on} \qquad B(R_i) \- B(R_0)\\
0& \qquad \text{on} \qquad M\- B(R_i)\endaligned \right.
$$
where $f_i$ is given by the sequence of harmonic functions obtained from \thetag{1.2} and \thetag{1.3}.  
Setting $\phi =\phi_i$ in the weighted Poincar\'e inequality, we have
$$
\split
\int_M \rho\,\phi_i^2 &\le \int_M |\n \phi_i|^2\\
&= \int_{B(R_i)\- B(R_0)} |\n f_i|^2\\
& = \int_{\p B(R_i)} f_i\, \frac{\p f_i}{\p \nu} - \int_{\p B(R_0)} f_i\, \frac{\p f_i}{\p \nu}\\
& = - \int_{\p B(R_0)} \frac{\p f_i}{\p \nu}.
\endsplit
$$
However, since $f_i \to 1$ on $M\- B(R_0)$, we conclude that 
$$
\int_M \rho \le 0,
$$
violating the assumption that $\rho \ge 0$ is nonzero.  This proves the second part of the corollary.\qed
\enddemo

While the existence of a weighted Poincar\'e inequality is equivalent to nonparabolicity, the condition that $M$ has property ($\Cal P_\rho$) is not as clear cut.  
The following lemma gives a sufficient condition for ($\Cal P_\rho$).

\proclaim{Lemma 1.5} Let $M$ be a nonparabolic complete manifold.  Suppose $G(x, y)$ is the minimal, symmetric, positive, Green's function for 
the Laplacian acting on $L^2$ functions.   For a fixed point $p \in M$, if $G(p, x) \to 0$ as $x \to \infty$, then $M$ has property ($\Cal P_\rho$) with
$$
\rho = \frac{|\n G(p, x)|^2}{4G^2(p, x)}.
$$
\endproclaim

\demo{Proof} In view of Corollary 1.4, it remains to show that the $\rho$-metric is complete.  Indeed, let $\gamma (s)$ be a curve parametrized by arclength,
$0\le s\le T.$  
The length of $\gamma$ with respect to the $\rho$-metric is given by
$$
\int_\gamma \sqrt{\rho}\, ds =\frac 12 \int_\gamma |\n \log G|\, ds.
$$
However, since
$$
\split
\log G(p, \gamma(0)) - \log G(p, \gamma(T)) &= -\int_0^T \frac{\p}{\p s}\left( \log G(p, \gamma(s)\right)) \, ds\\
& \le \int_\gamma |\n \log G|\,ds,
\endsplit
$$
we conclude that the $\rho$-length of $\gamma$ is infinity if $\gamma(T) \to \infty.$  This proves the completeness of the $\rho$-metric. \qed
\enddemo

We should point out that it is not necessarily true that if $G(p, x)$ does not tend to 0 at infinity then $M$ does not have property ($\Cal P_\rho$) 
since there might be other weight function that gives a complete metric.

\noindent
{\bf Example 1.6.} When $M= \Bbb R^n$ for $n \ge 3,$ the Green's function is given by
$$
G(0, x) = C_n\,r^{2-n}(x)
$$
for some constant $C_n >0$ depending only on $n.$  In this case, we compute that
$$
\frac{|\n G(0, x)|^2}{4 G^2(0, x)} = \frac{(n-2)^2}{4}\, r^{-2}(x).
$$
Hence from the above discussion, we conclude that the weighted Poincar\'e inequality
$$
\frac{(n-2)^2}{4}\int_{\Bbb R^n} r^{-2}\, \phi^2 \le \int_{\Bbb R^n} |\n \phi|^2,
$$
which is the well-known Hardy's inequality, must be valid for all compactly supported smooth function  
$\phi \in C_c^{\infty}(\Bbb R^n)$, and $\Bbb R^n$ has property ($\Cal P_\rho$) with
$$
\rho (x) = \frac{(n-2)^2}4\,r^{-2}(x).
$$

\noindent
{\bf Example 1.7.} Let $M^n$ be a minimal submanifold of dimension $n \ge 3$  in $\Bbb R^N.$  If we denote $\bar r$ to be the extrinsic 
distance function of $\Bbb R^N$ to a fixed point $p \in M,$ then it is known that it satisfies the equation
$$
\d \bar r \ge (n-1)\,\bar r^{-1},\tag 1.4
$$
where $\d$ is the Laplacian on $M$ with respect to the induced metric from $\Bbb R^N.$
For any $\phi \in C_c^{\infty}(M)$, we consider the integral
$$
\split
\int_{M} \bar r^{-1}\, \phi^2\, \d \bar r & = -2\int_M \bar r^{-1}\, \phi \la \n \phi, \n \bar r \ra + \int_M \bar r^{-2} \phi^2 |\n \bar r|^2\\
&\le 2 \int_M \bar r^{-1}\, \phi\, |\n \phi| + \int_M \bar r^{-2} \phi^2.
\endsplit
$$
Combining with \thetag{1.4}, this implies that
$$
\split
\frac{(n-2)}2 \int_M \bar r^{-2}\, \phi^2 &\le \int_M \bar r^{-1}\, \phi\, |\n \phi|\\
& \le \left( \int_M \bar r^{-2}\, \phi^2 \right)^{\frac 12} \left( \int_M |\n \phi|^2 \right)^{\frac 12},
\endsplit
$$
hence the weighted Poincar\'e inequality is valid on $M$ with 
$$
\rho = \frac {(n-2)^2}4 \, \bar r^{-2}.
$$
If we further assume that $M$ is properly immersed, then Lemma 1.5 implies that $M$ has property ($\Cal P_\rho$).

\noindent
{\bf Example 1.8.} Let $M$ be a simply connected, complete, Cartan-Hadamard manifold with sectional curvature bounded from above by 
$$
K_M \le -1.
$$
In this case, the Hessian comparison theorem asserts that
$$
\d r \ge (n-1) \coth r, \tag 1.5
$$
where $r$ is the geodesic distance function from a fixed point $ p \in M.$  Suppose $\phi \in C_c^{\infty}(M)$ is a compactly supported function.
Then \thetag{1.5} and integration by parts yield
$$
\split
(n-1) \int_M \phi^2\, \coth r &\le \int_M \phi^2\, \d r\\
&= -2\int_M \phi\, \la \n \phi, \n r\ra\\
&\le 2\int_M \phi\, |\n \phi|\\
& \le \frac{n-1}2 \int_M \phi^2 + \frac{2}{n-1} \int_M |\n \phi|^2.
\endsplit
$$
This can be rewritten as
$$
\frac{(n-1)^2}4 \int_M \phi^2 + \frac{(n-1)^2}2 \int_M (\coth r -1)\, \phi^2 \le \int_M |\n \phi|^2.
$$
Hence in this case, the weight function $\rho$ can be taken to be
$$
\split
\rho &= \frac{(n-1)^2}4 + \frac{(n-1)^2}{2} \, (\coth r -1)\\
&\ge \frac{(n-1)^2}{4}.
\endsplit
$$
Since it is bounded from below by a positive constant, the $\rho$-metric must be complete by the completeness assumption on $ds^2_M,$ and $M$ has property ($\Cal P_\rho$).

\heading
\S2 Decay Estimate
\endheading

In this section, we will consider a more general situation.
Let $V$ be a given potential function defined on $M,$ and 
$$
\d -V(x)
$$
be the Schr\"odinger operator on $M.$  We assume that there exists a positive function $\rho$ defined on $M$, such that, the weighted Poincar\'e type inequality
$$
\int_M \rho(x)\, \phi^2(x)\,dx \le \int_M |\n \phi|^2(x)\,dx + \int_M V(x)\,\phi^2(x)\,dx\tag 2.1
$$
is valid for any compactly supported function $\phi \in C_c^{\infty}(M).$ 
Let us define the $\rho$-metric  given by
$$
ds_{\rho}^2 =\rho\, ds^2.
$$
Using this metric, we consider the $\rho$-distance function defined to be,
$$
r_\rho (x, y) = \inf_{\gamma} \ell_{\rho}(\gamma),
$$
the infimum of the lengths of all smooth curves joining $x$ and $y$ with respect to $ds_{\rho}^2.$  
For a fixed point $p\in M$, we denote $r_\rho (x) = r_\rho(p, x)$ to be the $\rho$-distance to $p.$  
One checks readily that $|\n r_\rho |^2(x) = \rho(x).$ As in the case when $V = 0,$ we say that the manifold has 
property ($\Cal P_{\rho,V}$) if the $\rho$-metric is complete, and this will be the standing assumption on that $M$.

Throughout this article, we denote 

$$
B_\rho(p, R)=\{ x \in M\,|\, r_\rho(p,x) < R\}
$$
to be the set of points in $M$ that has $\rho$-distance less than $R$ from point $p.$    We also denote

$$
B(p, R) = \{ x \in M\,|\, r(p,x) < R\}
$$
to be the set of points in $M$ that has distance less than $R$ from point $p$ with respect to the background metric $ds_M^2.$  
When $p\in M$ is a fixed point, we will suppress the dependency of $p$ and write
$B_\rho(R) = B_\rho(p, R)$ and $B(R) = B(p, R).$
If $E$ is an end of $M$, we denote $E_\rho(R) = B_\rho(R) \cap E.$   

\proclaim{Theorem 2.1} Let $M$ be a complete Riemannian manifold with property ($\Cal P_{\rho,V}$).
Let $f$ be a nonnegative function defined on $E$ satisfying the differential inequality
$$
(\d - V(x)) \,f(x)\ge 0.
$$
If $f$ satisfies the growth condition
$$
\int_{E_\rho(R)} \rho\, f^2\, \exp(-2r_\rho ) = o(R)
$$
as $R \to \infty$, then it
must satisfy the decay estimate
$$
\int_{E_\rho(R+1)\setminus E_\rho(R)} \rho\, f^2 \le C\,\exp(-2R)
$$
for some constant $C>0$ depending on $f$ and $\rho$.
\endproclaim

\demo{Proof} We will first prove that for any $0<\delta<1,$ there exists a constant $0< C< \infty$ such that,
$$
\int_{E} \rho\,\exp(2\delta r_\rho)\, f^2 \le C.
$$
Indeed, let $\phi(r_\rho (x))$ be a non-negative cut-off function with support in $E$
with $r_\rho (x)$ being the $\rho$-distance to the fixed point $p.$  Then for any function $h(r_\rho (x))$ integration by parts yields
$$ 
\split 
\int_{E} &|\n (\phi \,f\,\exp (h))|^2\\
&\quad = \int_{E} |\n (\phi\,\exp(h))|^2\,f^2 + \int_{E} (\phi\,\exp(h))^2\,|\n f|^2 \\
& \quad \qquad +2\int_{E} \phi\,f\,\exp(h)\,\la \n(\phi\,\exp(h)), \n f \ra\\ 
& \quad = \int_{E} |\n (\phi\,\exp(h))|^2\,f^2   + \int_{E} \phi^2\,\exp(2h)\,|\n f|^2 \\
&\quad \qquad + \frac 12 \int_{E} \la \n(\phi^2\,\exp(2h)), \n (f^2) \ra\\ 
& \quad = \int_{E} |\n (\phi\,\exp(h))|^2\,f^2  + \int_{E} \phi^2\,\exp(2h)\,|\n f|^2  \\
&\quad \qquad - \frac 12 \int_{E} \phi^2\,\exp(2h)\, \d (f^2)\\ 
& \quad =\int_{E} |\n (\phi\,\exp(h))|^2\,f^2 - \int_{E} \phi^2\,\exp(2h)\, f\,\d f\\
& \quad \le \int_{E} |\n (\phi\,\exp(h))|^2\,f^2 - \int_{E} V\,\phi^2\,\exp(2h)\, f^2\\
&  \quad =\int_{E} |\n \phi|^2\,f^2\, \exp(2h) + 2 \int_{E} \phi\, \exp(2h)\, \la \n \phi, \n h\ra f^2\\
& \quad \qquad  +  \int_{E} \phi^2\,|\n h|^2\,  f^2\, \exp(2h) - \int_{E} V\,\phi^2\, f^2\,\exp(2h).
\endsplit \tag2.2
$$
On the other hand, using the  assumption \thetag{2.1}, we have
$$ 
\int_{E} \rho\,\phi^2 \,f^2 \exp(2h) \le \int_{E} |\n
(\phi \,f\, \exp(h))|^2 + \int_E V\,\phi^2\, f^2 \exp(2h), 
$$ 
hence \thetag {2.2} becomes
$$
\split
\int_{E}& \rho\, \phi^2\, f^2\, \exp(2h)\\
& \le   \int_{E} |\n \phi|^2\, f^2\, \exp(2h) + 2 \int_{E} \phi\,\exp(2h) \, \la \n \phi, \n h\ra f^2  +  \int_{E} \phi^2\,|\n h|^2\, f^2\, \exp(2h).
\endsplit\tag 2.3
$$

Let us now choose
$$
\phi (r_\rho (x))= \left\{ \aligned r_\rho (x)-R_0& \qquad \text{on} \qquad E_\rho(R_0+1) \- E_\rho(R_0)\\
1\qquad& \qquad \text{on} \qquad E_\rho(R)\setminus E_\rho(R_0+1),\\
R^{-1}(2R -&r_\rho (x)) \qquad \text{on } \qquad E_\rho(2R)\- E_\rho(R)\\
0 \qquad & \qquad \text{on } \qquad E\- E_\rho(2R),
\endaligned \right.
$$
and hence
$$
|\n \phi|^2 (x)= \left\{ \aligned \rho(x)\quad & \qquad \text{on} \qquad E_\rho(R_0+1) \- E_\rho(R_0)\\
R^{-2}\,\rho(x)& \qquad \text{on } \qquad E_\rho(2R)\- E_\rho(R)\\
0 \qquad & \qquad \text{on } \qquad (E_\rho(R) \- E_\rho(R_0 +1)) \cup (E\- E_\rho(2R)).
\endaligned \right.
$$
We also choose
$$
h(r_\rho (x)) = \left\{ \aligned \delta \,& r_\rho (x) \qquad \text{for } \qquad r_\rho \le \frac{K}{(1+\delta) } \\
K& - r_\rho (x) \qquad \text{for } \qquad r_\rho \ge \frac{K}{(1+\delta)} 
\endaligned \right.
$$
for  some  fixed $K >(R_0+1)(1+ \delta).$
When $R \ge \frac{K}{(1+\delta)}$, we see that
$$
|\n h|^2(x) = \left\{ \aligned \delta^2& \rho(x) \qquad \text{for } \qquad r_\rho \le \frac{K}{(1+\delta)}\\
\rho&(x) \qquad \text{for } \qquad r_\rho \ge \frac{K}{(1+\delta) } 
\endaligned \right.
$$
and
$$
\la \n \phi, \n h \ra (x) = \left\{ \aligned \delta\, \rho&(x) \qquad \text{on } \qquad E_\rho(R_0+1)\-E_\rho(R_0)\\
R^{-1} &\rho (x) \qquad \text{on } \qquad E_\rho(2R) \-E_\rho(R)\\
 0& \qquad \text{otherwise.}\endaligned \right.
$$
Substituting into \thetag{2.3}, we obtain
$$
\split
\int_{E} &\rho\,\phi^2\,  f^2\, \exp(2h)\\
& \le   \int_{E_\rho(R_0+1) \- E_\rho(R_0)}  \rho\,f^2\, \exp(2h) + R^{-2}\int_{E_\rho(2R) \- E_\rho(R)} \rho\, f^2 \, \exp(2h)\\
&\qquad + 2\delta  \int_{E_\rho(R_0+1) \- E_\rho(R_0)} \rho \, f^2\, \exp(2h)\\
& \qquad +2R^{-1} \int_{E_\rho(2R) \- E_\rho(R)} \rho\, f^2\, \exp(2h)\\
&\qquad + \delta^2 \int_{E_\rho(K(1+\delta)^{-1}) \- E_\rho(R_0))} \rho\, \phi^2\,f^2\, \exp(2h)\\
&  \qquad  +  \int_{E_\rho(2R) \- E_\rho(K(1+\delta)^{-1})} \rho\, \phi^2\,f^2 \, \exp(2h).
\endsplit
$$
This can be rewritten as
$$
\split
\int_{E_\rho(K(1+\delta)^{-1}) \- E_\rho(R_0+1)} &\rho\, f^2\, \exp(2h)\\
& \le \int_{E_\rho(K(1+\delta)^{-1})} \rho\,\phi^2\,  f^2\, \exp(2h) \\
& \le \int_{E_\rho(R_0+1) \- E_\rho(R_0)}  \rho\,f^2\, \exp(2h)\\ 
& \qquad + R^{-2}\int_{E_\rho(2R) \- E_\rho(R)} \rho\, f^2\, \exp(2h)\\
&\qquad + 2\delta  \int_{E_\rho(R_0+1) \- E_\rho(R_0)} \rho\, f^2\, \exp(2h)\\
& \qquad +2R^{-1}  \int_{E_\rho(2R) \- E_\rho(R)} \rho\, f^2\, \exp(2h)\\
&\qquad + \delta^2 \int_{E_\rho(K(1+\delta)^{-1}) \- E_\rho(R_0))} \rho \,f^2\, \exp(2h),
\endsplit
$$
hence
$$
\split
(1- \delta^2)&\int_{E_\rho (K(1+\delta)^{-1}) \- E_\rho(R_0+1))}  \rho\, f^2\, \exp(2h)\\
& \le   ( \delta^2 + 2 \delta + 1 ) \int_{E_\rho(R_0+1) \- E_\rho(R_0)} \rho\,f^2\, \exp(2h) \\
& \qquad  +R^{-2} \int_{E_\rho(2R) \- E_\rho(R)} \rho\, f^2 \, \exp(2h) + 2R^{-1}  \int_{E_\rho(2R) \- E_\rho(R)} \rho\,  f^2\, \exp(2h).
\endsplit
$$
The definition of $h$ and the assumption on the growth condition on $f$ imply that the last two terms on the right hand side tend 
to $0$ as $ R \to \infty.$  Hence we obtain the estimate
$$
\split
(1- \delta^2) &\int_{E_\rho(K(1+\delta)^{-1}) \- E_\rho(R_0+1))} \rho\, f^2\, \exp(2\delta r_\rho)\\
& \le   ( \delta^2 + 2 \delta + 1) \int_{E_\rho(R_0+1) \- E_\rho(R_0))} \rho\,f^2\, \exp(2\delta r_\rho).
\endsplit
$$
Since the right hand side is independent of $K$, by letting $ K \to \infty$ we conclude that
$$
\int_{E \-E_\rho(R_0+1)} \rho \,f^2\,\exp(2\delta  r_\rho)  \le C, \tag 2.4
$$
for some constant $0 < C < \infty.$

Our next step is to improve this estimate by setting $h= r_\rho$
in the preceding argument. Note that with this choice of $h$,  \thetag{2.3} asserts that
$$
-2\,\int_E \phi\,\exp(2r_\rho) \,\la \n \phi, \n r_\rho\ra \,f^2 \le
\int_E |\n \phi|^2\,f^2\,\exp(2r_\rho).
$$
For $R_0< R_1< R$, let us choose $\phi$ to be
$$ \phi(x) = \left\{ \aligned  \frac{r_\rho(x) - R_0}{R_1 - R_0} \qquad
&\text{on} \qquad E_\rho(R_1)\setminus E_\rho(R_0)\\ \frac{R - r_\rho(x)}{R -
R_1} \qquad & \text{on} \qquad E_\rho(R) \setminus E_\rho(R_1).
\endaligned \right.
$$ 
We conclude that
$$
\split
 \frac{2}{R-R_1}&\,\int_{E_\rho(R) \- E_\rho(R_1)}
\left(\frac{R - r_\rho(x)}{R - R_1}\right)\,\rho\,f^2\,\exp(2r_\rho) \\
&\le \frac{1}{(R_1 - R_0)^2} \int_{E_\rho(R_1)\setminus E_\rho(R_0)} \rho\,f^2 \, \exp(2r_\rho) \\
& \qquad + \frac{1}{(R-R_1)^2}\,\int_{E_\rho(R) \- E_\rho(R_1)} \rho \,f^2\,\exp(2r_\rho) \\
&\qquad + \frac{2}{(R_1 - R_0)^2} \int_{E_\rho(R_1)\- E_\rho(R_0)} (r_\rho-
R_0)\,\rho \,f^2\,\exp(2r_\rho) .
\endsplit
$$
On the other hand, for any $0 < t<R -R_1$, since
$$
\split
& \frac{2t}{(R-R_1)^2}\,\int_{E_\rho(R-t) \setminus E_\rho(R_1)} \rho \,f^2\, \exp(2r_\rho) \\
& \qquad\le \frac{2}{(R-R_1)^2}\,\int_{E_\rho(R) \setminus E_\rho(R_1)}
(R - r_\rho(x))\, \rho \,f^2\, \exp(2r_\rho),
\endsplit
$$
we deduce that
$$
\split
&\frac{2t}{(R-R_1)^2}\,\int_{E_\rho(R-t) \setminus E_\rho(R_1)} \rho \,f^2\,\exp(2r_\rho) \\
& \qquad  \le \left(\frac{2}{R_1 - R_0} + \frac 1{(R_1-R_0)^2} \right) \int_{E_\rho(R_1)\setminus
E_\rho(R_0)} \rho\,f^2\, \exp(2r_\rho) \\& \qquad \quad +
\frac{1}{(R-R_1)^2}\int_{E_\rho(R)
\- E_\rho(R_1)} \rho \,f^2\,\exp(2r_\rho).
\endsplit \tag 2.5
$$ 

Observe that by taking  $R_1 = R_0+1$, $t =1$, and setting
$$
g(R) = \int_{E_\rho(R) \setminus E_\rho(R_0+1)} \rho\,f^2\,\exp(2r_\rho),
$$
the inequality \thetag {2.5} can be written as
$$
g(R-1) \le C_1\,R^2 + \frac 12\, g(R),
$$
where
$$
C_1 = \frac{3}{2} \int_{E_\rho(R_0 + 1) \setminus E_\rho(R_0)} \rho\,f^2\,\exp(2r_\rho)
$$
is independent of $R.$ Iterating this inequality, we obtain that for
any positive integer $k$ and $R\ge 1$
$$
\split
g(R) &\le C_1\,\sum_{i=1}^{k}\frac{(R+i)^2}{2^{i-1}} + 2^{-k}\,g(R+k)\\
&\le C_1\,R^2 \sum_{i=1}^{\infty}\frac{(1+i)^2}{2^{i-1}} + 2^{-k}\,g(R+k)\\
&\le C_2\,R^2 + 2^{-k}\,g(R+k)
\endsplit
$$ 
for some constant $C_2.$
However, our previous estimate \thetag {2.4}
asserts that
$$
\int_{E} \rho \,f^2\,\exp(2\delta r_\rho) \le C
$$
for any $\delta < 1.$  This implies that
$$ 
\split 
g(R+k) &= \int_{E_\rho(R+k) \setminus E_\rho(R_0+1 )} \rho\,f^2\, \exp(2r_\rho)\\
 &\le \exp(2(R+k)(1-\delta)) \int_{E_\rho(R+k) \setminus
E_\rho(R_0+1 )} \rho\,f^2\,\exp(2\delta r_\rho)\\ &\le
C\,\exp(2(R+k)(1 - \delta)).
\endsplit
$$
Hence,
$$
2^{-k}\,g(R+k) \to 0
$$
as $k \to \infty$ by choosing $2(1- \delta) < \ln 2.$
This proves the estimate 

$$
g(R)\le C_2\, R^2.
$$
By adjusting the constant, we have

$$
\int_{E_\rho(R)} \rho\,f^2\,\exp(2r_\rho) \le C_3\,R^2 \tag 2.6
$$
for all $R\ge R_0.$

Using inequality \thetag {2.5} again and by  choosing $R_1 = R_0 +1$ and $t = \frac R2$ this time, we conclude that
$$
\split
R&\,\int_{E_\rho(\frac R2) \setminus E_\rho(R_0+1)} \rho\,f^2\,\exp(2r_\rho) \\
& \qquad  \le C_4\,R^2 + \int_{E_\rho(R) \setminus E_\rho(R_0+1)} \rho\,f^2\,\exp(2r_\rho).
\endsplit
$$
However, applying the estimate \thetag{2.6} to the second term
on the right hand side, we have
$$
\int_{E_\rho(\frac R2) \setminus E_\rho(R_0+1)} \rho\,f^2\,
\exp(2r_\rho) \le C_5\, R.
$$
Therefore, for $R\ge R_0$,
$$
\int_{E_\rho(R)} \rho \,f^2\,\exp(2r_\rho) \le C\, R. \tag 2.7
$$

We are now ready to prove the theorem by using \thetag {2.7}.
Setting $t= 2$ and $R_1 = R - 4$ in \thetag{2.5}, we obtain

$$
\split
&\int_{E_\rho(R-2) \setminus E_\rho(R-4)} \rho\,f^2 \,\exp(2r_\rho) \\
& \qquad  \le \left(\frac{8}{R - R_0 -4} + \frac 4{(R - R_0 -
4)^2}\right)\,\int_{E_\rho(R-4)\setminus
E_\rho(R_0)} \rho \,f^2\,\exp(2r_\rho) \\
& \qquad \quad + \frac{1}{4}\int_{E_\rho(R)\setminus E_\rho(R-4)} \rho \,f^2\, \exp(2r_\rho).
\endsplit
$$
According to \thetag {2.7}, the first term of the right hand side is bounded by a 
constant.
Hence, the above inequality can be rewritten as
$$ 
\int_{E_\rho(R-2) \setminus E_\rho(R-4)} \rho \,f^2\, \exp(2r_\rho) \le C +
\frac{1}{3}\int_{E_\rho(R) \setminus E_\rho(R-2)} \rho \,f^2\,\exp(2r_\rho). 
$$
Iterating this inequality $k$ times, we arrive at
$$
\split
\int_{E_\rho(R+2) \setminus E_\rho(R)} &\rho \,f^2\,\exp(2r_\rho)\\ 
& \le C \,\sum_{i=0}^{k-1} 3^{-i} + 3^{-k} \int_{E_\rho(R + 2(k+1))\setminus E_\rho(R + 2k)} \rho \,f^2\, \exp(2r_\rho).
\endsplit
$$
However, using \thetag{2.7} again, we conclude that the second term is bounded by
$$
3^{-k} \int_{E_\rho(R + 2(k+1))\setminus E_\rho(R + 2k)} \rho \,f^2\,\exp(2r_\rho) \le C\,3^{-k} (R + 
2(k+1))
$$
which tends to $0$ as $k\to \infty.$  Hence
$$
\int_{E_\rho(R + 2) \setminus E_\rho(R)} \rho \,f^2\, \exp(2r_\rho) \le C. \tag 2.8
$$
for some constant $C>0$ independent of $R.$ The theorem now follows from \thetag 
{2.8}. \qed
\enddemo

We now draw some corollaries. 

\proclaim{Corollary 2.2} Let $M$ be a complete Riemannian manifold.
Suppose $E$ is an end of $M$ such that $\lambda_1(E) >0$, i.e., 
$$
\lambda_1(E)\,\int_E \phi^2(x)\,dx \le \int_E |\n \phi|^2(x)\,dx
$$
for any compactly supported function $\phi \in C_c^{\infty}(E).$  Let $f$ be a nonnegative function defined on $E$ satisfying the differential inequality
$$
(\d + \mu) \,f(x)\ge 0
$$
for some constant $\mu$ with the property that $\lambda_1(E) - \mu >0.$
If $a = \sqrt{\lambda_1(E) - \mu}$ and $f$ satisfies the growth condition
$$
\int_{E(R)}  f^2\, \exp(-2ar) = o(R)
$$
as $R \to \infty$, then it
must satisfy the decay estimate
$$
\int_{E(R+1)\setminus E(R)} f^2 \le C\,\exp(-2aR)
$$
for some constant $C>0$ depending on $f$ and $a$. 
\endproclaim

\demo{Proof} By setting $-V(x) = \mu$ we can rewrite the Poincar\'e inequality as
$$
(\lambda_1(E) - \mu)  \int_E \phi^2(x)\,dx \le \int_E |\n \phi|^2(x)\, dx - \mu \int_E \phi^2(x)\,dx.
$$
We now can apply Theorem 2.1 by setting $\rho = a^2.$  The distance function with respect to the metric $\rho\, ds^2$ is then given by
$$
r_\rho (x)= a\, r(x)
$$
where $r(x)$ is the background distance function to the smooth compact set $\Omega \subset M.$    The corollary follows from Theorem 2.1.\qed
\enddemo

\proclaim{Corollary 2.3} Let $M$ be a complete Riemannian manifold satisfying property ($\Cal P_\rho$).  Suppose $\{E_1, \dots, E_k\}$ with $k \ge 2$ 
are the nonparabolic ends of $M$. Then for each $1\le i \le k$ there exists a bounded harmonic function $f_i$ defined on $M$ satisfying the growth estimate
$$
\int_{B_\rho(R+1) \- B_\rho(R)} |\n f_i|^2 \le C\, \exp(-2R).
$$ 
Moreover, $0 \le f_i \le 1$ and has the property that
$$
\sup_{x \in E_i} f_i (x) = 1,
$$
and
$$
\inf_{x \in E_j} f_i(x) = 0, \qquad \text{for} \qquad  j \neq i.
$$
\endproclaim

\demo{Proof}We will construct $f_i$ for the case $i=1$, and  the construction for other values of $i$ is exactly the same.  
In this case, we will simply denote $f = f_1.$  Following the theory of Li-Tam \cite{L-T2} (see also \cite{L-W1}), $f$ can be constructed by 
taking the limit, as $\bar R \to \infty$, of a converging subsequence of harmonic functions $f_{\bar R}$ satisfying 
$$
\d f_{\bar R} = 0 \qquad \text{on} \qquad B(\bar R),
$$
with boundary condition
$$
f_{\bar R} = 1 \qquad \text{on} \qquad \p B(\bar R) \cap E_1,
$$
and
$$
f_{\bar R} = 0  \qquad \text{on} \qquad \p B(\bar R) \- E_1.
$$
In fact, we only need to verify the growth estimate for the Dirichlet integral for the limiting function.  The other  required 
properties of $f$ follow from the contruction of Li-Tam.  To check the growth estimate, we first show that 
on $(B_\rho(R+1) \- B_\rho(R)) \-E_1$, because of the boundary condition  we can apply Theorem 2.1 to the function $f_{\bar R}.$  
By taking the limit, this implies that
$$
\int_{(B_\rho(R+1) \- B_\rho (R)) \- E_1} \rho\,f^2 \le C\, \exp(-2R).\tag 2.9
$$
Similarly on $E_1,$ we can apply Theorem 2.1 to the function $1 - f_{\bar R},$ hence we obtain
$$
\int_{(E_1 \cap B_\rho(R+1)) \- (E_1 \cap B_\rho(R))} \rho\, (1-f)^2 \le C\, \exp(-2R).\tag 2.10
$$
Let us now consider the cut-off function
$$
\phi (r_{\rho}(x)) = \left\{ \aligned r_\rho(x)& - R+1 \qquad \text{for} \qquad R-1 \le r_\rho \le R\\
&1 \qquad \text{for} \qquad R \le r_\rho \le R+1\\
R+2 &- r_\rho \qquad \text{for} \qquad R+1 \le r_\rho \le R+2\\
0& \qquad \text{otherwise}.
\endaligned \right.
$$
Integrating by parts and Schwarz inequality yield
$$
\split
0 & = \int_{(B_\rho(R+2) \- B_\rho (R-1)) \- E_1} \phi^2 \, f \d f\\
&= - \int_{(B_\rho(R+2) \- B_\rho (R-1)) \- E_1} \phi^2 \, |\n f|^2 - 2 \int_{(B_\rho(R+2) \- B_\rho (R-1)) \- E_1} \phi\,f\,\la \n \phi, \n f \ra\\
&\le - \frac 12\int_{(B_\rho(R+2) \- B_\rho (R-1)) \- E_1} \phi^2 \, |\n f|^2 + 2 \int_{(B_\rho(R+2) \- B_\rho (R-1)) \- E_1} |\n \phi|^2\,f^2.
\endsplit
$$
Hence combining with the definition of $\phi$, we obtain the estimate
$$
\split
 \int_{(B_\rho(R+1) \- B_\rho (R)) \- E_1}  |\n f|^2 &\le  \int_{(B_\rho(R+2) \- B_\rho (R-1)) \- E_1} \phi^2 \, |\n f|^2 \\
& \le 4 \int_{(B_\rho(R+2) \- B_\rho (R-1)) \- E_1} \rho\, f^2.
\endsplit
$$
Applying the estimate \thetag {2.9} to the right hand side, we conclude the desired estimate on the set $(B_\rho(R+1) \- B_\rho (R)) \- E_1$.  
The estimate on $(E_1 \cap B_\rho(R+1)) \-(E_1 \cap B_\rho(R))$ follows by using the function $1-f$ and \thetag{2.10} instead. \qed

\enddemo

We would like to point out  that the hypothesis of Corollary 2.2, hence Theorem 2.1, is best possible.  
Indeed, if we consider the hyperbolic space form $\Bbb H^n$ of $-1$ constant sectional curvature, then the volume growth is given by
$$
V(R) \sim C\, \exp((n-1)R)
$$
and 
$$
\lambda_1(\Bbb H^n) = \frac{(n-1)^2}{4}.
$$
We consider Theorem 2.1 for the special case when $V(x)=0$ and $\rho = \frac{(n-1)^2}4.$ In this case, the distance function $r_\rho$ 
with respect to the metric $\rho \,ds^2$ is simply given by
$$
r_\rho = \frac{(n-1) }2\, r,
$$
where $r$ is the hyperbolic distance function.
If $f$ is a nonconstant bounded harmonic function, then
$$
\int_{B(R)} \rho\, f^2\,\exp(-2r_\rho) = O(R).
$$
We claim that the conclusion of Theorem 2.1 is not valid, hence will imply that the hypothesis of Theorem 2.1 cannot be improved.  
Indeed, if the conclusion were true, then $f$ would be in $L^2(\Bbb H^n)$.  However, Yau's theorem \cite{Y2} implies that $f$ must be identically constant.  
On the other hand, it is known that $\Bbb H^n$ has an infinite dimensional space of bounded harmonic functions, which provides a contradiction.

Also note that in the case of $\Bbb R^n$ $(n \ge 3)$,  the distance function with respect to the $\rho$-metric is
$$
\split
r_\rho &\sim \frac{n-2}2 \int_1^r t^{-1}\,dt\\
&= \log r^{\frac{n-2}2}.
\endsplit
$$
as $ r \to \infty.$  If we consider a multiple of the Green's function $f( x) = r^{2-n}$ on $\Bbb R^n$, then checking the 
hypothesis of Theorem 2.1 for $f(x)$ on $E = \Bbb R^n \- B(1),$ the integral
$$
\split
\int_E \rho\,\,f^2\, \exp(-2r_\rho) &=\frac{(n-2)^2}4 \int_1^{\infty} r^{-2}\,r^{-n+2}\,r^{4-2n}\,r^{n-1}\, dr\\
& = \frac{(n-2)^2}4 \int_1^{\infty} r^{-2n +3}\,dr\\
& < \infty.
\endsplit
$$
Hence we can apply Theorem 2.1 to this choice of $f$.  On the other hand, the integral
$$
\split
\int_{E_\rho (R+1)\- E_\rho (R)} \rho\, f^2 &=\frac{(n-2)^2}4 \int_{r_\rho = R+1}^{r_\rho =R} r^{-2}\,r^{4-2n}\,r^{n-1}\, dr\\
& = \int_{e^{\frac{2(R+1)}{n-2}}}^{e^{\frac{2R}{n-2}}} r^{-n +1}\,dr\\
& = \frac{ 1 }{n-2} (\exp(-2R) - \exp(-2(R+1)))\\
& = \frac{1- e^{-2}}{n-2}\, \exp(-2R).
\endsplit
$$
This implies that the conclusion of Theorem 2.1 is also sharp in this case.

Finally, we point out that the preceding argument of Theorem 2.1 can be extended without much modification 
to deal with $p$-forms satisfying a suitable differential equation. We consider the operator
$$
\d +W(x)
$$
acting on the $p$-forms on $M,$ where $\d$ is the Hodge Laplacian and $W$ an endomorphism on 
the bundle of $p$-forms on $M.$ 

\proclaim{Theorem 2.4} Let $M$ be a complete Riemannian manifold.
Suppose $E$ is an end of $M$ such that there exists a nonnegative function $\rho(x)$ defined on $E$ with the property that

$$
\int_E \rho(x)\, |\eta|^2(x)\,dx \le \int_E \left(|d \eta|^2(x)+|\delta \eta|^2(x)\right)\,dx + \int_E \la W(\eta)(x),\eta(x) \ra\,dx
$$
is valid for any compactly supported smooth $p$-form  $\eta$ on $E.$  Assume that the $\rho$-metric given by $ds_\rho^2 = \rho\, ds_M^2$ is complete on $E.$
Let $\omega$ be a smooth $p$-form defined on $E$ satisfying the differential inequality
$$
\la (\d +W(x))\omega, \omega\ra (x)\le 0
$$
for all $x\in E.$
If $\omega$ satisfies the growth condition
$$
\int_{E_\rho(R)} \rho\, |\omega|^2\, \exp(-2r_\rho ) = o(R)
$$
as $R \to \infty$, then it
must satisfy the decay estimate
$$
\int_{E_\rho(R+1)\setminus E_\rho(R)} \rho\, |\omega|^2 \le C\,\exp(-2R)
$$
for some constant $C>0$ depending on $\omega$ and $\rho$. 
\endproclaim

\heading
\S3 Geometric Conditions For Parabolicity and Nonparabolicity
\endheading

In this section, we would like to discuss some geometric conditions for the parabolicity and nonparabolicity of an end $E.$   In \cite{L-W1}, 
we used the decay estimate similar to \S2 to derive geometric conditions for parabolicity on a manifold with $\lambda_1(M) >0.$  
A similar argument will yield the following conditions for manifolds with property ($\Cal P_\rho$).  The key issue is that when ($\Cal P_\rho$) is present, 
the geometric conditions involving $\rho$ for parabolicity and nonparabolicity has a substantial gap.  This fact is important to the proof of our 
main theorems in the proceeding sections.

\proclaim{Theorem 3.1}Let $E$ be an end of a complete Riemannian manifold $M$ with property ($\Cal P_\rho$) for some weight functio $\rho.$  
If $E$ is nonparabolic, then
$$
\int_{E_\rho (R+1) \- E_\rho(R)} \rho\,dV \ge C_1\, \exp(2R)
$$
for some constant $C_1 >0$ and for $R$ sufficiently large, where $E_\rho(R) = B_\rho(R) \cap E.$ 
If $E$ is  parabolic, then
$$
\int_E \rho\,dV < \infty
$$
and
$$
\int_{E \- E_\rho(R)} \rho\, dV \le C_2\, \exp(-2R),
$$
for some constant $C_2>0$ and for all $R$ sufficiently large. 
\endproclaim

\demo{Proof} Following a similar argument as in the  proof of Theorem 1.4 in \cite{L-W1} and applying the harmonic equation to 
the barrier function $f$ on a nonparabolic end, we obtain
$$
\split
C  &= \int_{\p E} \frac{\p f}{\p \nu}\,dA\\
& = \int_{\p B_\rho(r) \cap E} \frac{\p f}{\p \nu}\,dA\\
& \le \left(\int_{\p B_\rho(r) \cap E} (\sqrt{\rho})^{-1}\, |\n f|^2 \,dA\right)^{\frac 12}\left( \int_{\p B_\rho(r) \cap E} \sqrt{\rho} \,dA \right)^{\frac 12}.
\endsplit \tag 3.1
$$
On the other hand, the co-area formula asserts that
$$
\split
\int_{E_\rho(R+1) \- E_\rho(R)} h\,dV &= \int_R^{R+1} \int_{\p B_\rho(r) \cap E} h\, |\n r_\rho|^{-1}\,dA\, dr_\rho\\
& = \int_R^{R+1} \int_{\p B_\rho(r) \cap E} h\, (\sqrt{\rho})^{-1} \,dA\,dr_\rho
\endsplit
$$
for any measurable function $h.$  Hence \thetag{3.1} together with Corollary 2.3 imply that
$$
\split
\int_R^{R+1} \left( \int_{\p B_\rho(r)\cap E} \sqrt{\rho} \, dA \right)^{-1} \, dr_\rho 
& \le C\, \int_R^{R+1} \int_{\p B_\rho (r) \cap E} |\n f|^2\, (\sqrt{\rho})^{-1} \, dA\, dr_\rho\\
&= C\, \int_{E_\rho(R+1) \- E_\rho(R)} |\n f|^2\, dV \\
& \le C\, \exp (-2R).
\endsplit
$$
Applying Schwarz inequality gives
$$
\split
1 &\le \left(\int_R^{R+1} \left( \int_{\p B_\rho(r) \cap E} \sqrt{\rho}\, dA \right)^{-1} dr_\rho \right)\,
\left(\int_R^{R+1} \int_{\p B_\rho(r) \cap E} \sqrt{\rho} \, dA\, dr_\rho \right) \\
& \le C\, \exp(-2R)\,\int_R^{R+1} \int_{\p B_\rho(r) \cap E} \sqrt{\rho} \, dA\, dr_\rho.
\endsplit \tag 3.2
$$
Using the co-area formula again, we obtain the estimate
$$
C_1\, \exp(2R) \le \int_{E_\rho (R+1) \- E_\rho(R)} \rho\,dV
$$
as claimed.

If $E$ is parabolic, we apply the proof of Corollary 2.3 to the barrier function $f = 1$ on E and obtain
$$
\int_{E_\rho(R+1) \- E_\rho(R)} \rho\, dV \le C\, \exp(-2R)
$$
for all sufficiently large $R.$
Summing over these estimate, we conclude that
$$
\split
\int_{E \- E_\rho(R)} \rho\, dV & = \sum_{i=0}^\infty \int_{E_\rho(R+i+1) \- E_\rho(R+i)} \rho\, dV\\
& \le C \sum_{i=0}^\infty \exp(-2(R+i))\\
& = C_2\, \exp(-2R)
\endsplit
$$
for some constant $C_2 >0.$  This proves the second half of the theorem.\qed
\enddemo

\proclaim{Corollary 3.2}Let $E$ be an end of a complete Riemannian manifold $M$ with property ($\Cal P_\rho$).  
If $E$ is nonparabolic then it must have at least quadratic volume growth.   In particular, if the weight function $\rho$  satisfies
$$
\liminf_{ x \to \infty} \rho(x) >0,
$$
then $E$ is nonparabolic if and only if $E$ has infinite volume.
\endproclaim

\demo{Proof} A Theorem of Varopoulos \cite{V} asserts that if $M$ is nonparabolic then
$$
\int_1^\infty A^{-1}(r)\,dr < \infty,
$$
where $A(r)$ denotes the area of the boundary of the geodesic ball of radius $r$ centered at a fixed point.  
In fact, this criterion can be localized \cite{L-T1} at an end, namely, that an end $E$ is nonparabolic implies that
$$
\int_1^\infty A_E^{-1}(r)\, dr < \infty, \tag 3.3
$$
where $A_E(r)$ denotes the area of the set $\p B(r) \cap E.$  In particular, applying the Schwarz inequality, we conclude that
$$
\split
R  &\le \left(\int_1^{R} A_E^{-1} (r)\, dr \right)^{\frac 12}\,\left( \int_1^{R} A_E(r)\, dr \right)^{\frac 12}\\
& \le C \,V_E^{\frac 12}(R),
\endsplit\tag 3.4
$$
where $V_E(R)$ is the volume of $E\cap B(R)$.
Hence $E$ must have infinite volume.
If $E$ is parabolic, then Theorem 3.1 implies that
$$
\int_E \rho\, dV < \infty.
$$
The assumption that $\liminf \rho >0$ implies that $E$ has finite volume. \qed
\enddemo

\proclaim{Corollary 3.3}Let $E$ be an end of a complete Riemannian manifold $M$ with property ($\Cal P_\rho$).  
If $\rho(r)$ is the weight function that depends only on the distance $r$ to a fixed compact set, then $E$ is nonparabolic if and only if
$$
\int_1^\infty A_E^{-1}(r)\, dr < \infty.
$$
\endproclaim

\demo{Proof} As pointed out in the proof of Corollary 3.2, \thetag{3.3} is a necessary condition for nonparabolicity.  
We now assume that $M$ is parabolic and by Theorem 3.1
$$
\int_E \rho\, dV < \infty.
$$
However, applying the Schwarz inequality, we have
$$
\split
r_\rho(R) &= \int_1^R \sqrt{\rho}(r)\, dr\\
&\le \left( \int_1^R \rho(r)\, A_E(r)\, dr \right)^{\frac 12} \left( \int_1^R A_E^{-1}(r)\, dr \right)^{\frac 12}.
\endsplit
$$
Letting $R \to \infty$ and using the completeness of the $\rho$-metric, we conclude that the left hand side tends to infinity, hence
$$
\infty = \int_1^\infty A_E^{-1}(r)\, dr,
$$
and the corollary is proved.\qed
\enddemo

\heading
\S4 Improved Bochner Formula and Metric Rigidity
\endheading

In this section, we will recall an improved Bochner formula and consider the case when the inequality is realized as an equality.  
This formula computes the Laplacian of the gradient of a harmonic function as in the standard Bochner formula, but with extra care 
taken on the Hessian term by using the harmonic equation one more time.  This manipulation was first used effectively by Yau \cite{Y1}, 
and it is by now considered to be a standard trick.

\proclaim{Lemma 4.1} Let $M^n$ be a complete Riemannian manifold of
dimension $n \ge 2.$   Assume that the Ricci curvature of $M$ satisfies the lower bound
$$
\text{Ric}_M(x) \ge - (n-1) \tau(x)
$$
for all $x \in M.$  Suppose $f$ is a nonconstant harmonic function defined on $M.$ Then the function  $|\n f|$ must satisfy the differential inequality
$$
\d |\n f| \ge -(n-1)\tau\, |\n f| + \frac{|\n |\n f| |^2}{(n-1) |\n f|}
$$ 
in the weak sense.  Moreover, if equality holds, then $M$ is given by $M = \Bbb R \times N^{n-1}$ with the warped product metric
$$
ds_M^2 = dt^2 + \eta^2(t)\,ds_N^2
$$
for some positive function $\eta(t)$, and some manifold $N^{n-1}.$ In this case, $\tau(t)$ is a function of $t$ alone satisfying
$$
\eta''(t)\, \eta^{-1}(t) = \tau(t).
$$
\endproclaim

\demo{Proof } If we denote 
$$
g = |\n f|,
$$ 
then the Bochner formula (see Theorem 2.1 of \cite{L-W1}) and the lower bound of the Ricci curvature assert that 
$$ 
\d g \ge -(n-1)\tau\, g + \frac{|\n g |^2}{(n-1) g}. \tag 4.1 
$$ 
Note that since $f$ is nonconstant, $g$ is not identically zero.  
Hence equality holds if and only if all the inequalities used in the proof of \thetag {4.1} are equalities.  
In particular, we conclude that there exists a function $\mu$, such that,
$$
f_{1\alpha} = 0\tag 4.2
$$
and
$$
f_{\alpha \beta} = \mu\, \delta_{\alpha \beta}\tag 4.3
$$  
for all $\alpha, \beta = 2, \dots ,n,$ where $\{e_1,e_2,\dots,e_n\}$ is an orthonormal frame satisfying
$|\n f|\,e_1=\n f$ and $e_\alpha\,f=0$ for all $\alpha\neq 1.$

We can now argue to conclude that $M = \Bbb R \times N$ with the warped product metric
$$
ds^2 = dt^2 + \eta^2(t) \,ds_N^2
$$
for some manifold $N$ and for some positive function $\eta$.

Indeed, since $\d f = 0$, together with \thetag{4.2} and \thetag{4.3}, the Hessian of $f$ must be
of the form
$$
(f_{ij}) = \left(\matrix -(n-1)\mu &0 &0 &\cdot  & \cdot &0\\
0 &\mu &0 &\cdot & \cdot &0\\
0 &0 &\mu &\cdot &\cdot &0\\
\cdot & \cdot & \cdot &\cdot\\
\cdot & \cdot &\cdot & &\cdot\\
0 &0 &0 &\cdot &\cdot &\mu
\endmatrix\right).
$$
The fact that $f_{1\alpha} = 0$ for all $\alpha \neq 1$ implies that $|\n f|$ is 
identically
constant along the level set of $f.$  In particular, the level sets of $|\n f|$ 
and $f$ coincide.  We claim that $|\n f|$ does not vanish anywhere hence $f$ has no critical points and $g >0.$  Assuming the contrary, if $|\n f|(x) = 0$, 
by addition of a constant, we may assume $f(x) = 0.$  The regularity theory of harmonic function asserts that $f$ locally in a neighborhood of $x$ 
behaves like a homogeneous harmonic polynomial in $\Bbb R^n$ with the origin at $x$.  This is impossible since the level set of $|\n f|$ and $f$ coincide.  
Hence $|\n f| >0$ and $M$ must be topologically the product $\Bbb R \times N$, where $N$ is given by the level set of $f.$ 
Also,
$$
\split
\mu\,\delta_{\alpha \beta} &= f_{\alpha \beta} \\
&= h_{\alpha \beta} \, f_1
\endsplit\tag 4.4
$$ with $(h_{\alpha \beta})$ being the second fundamental form of
the level set of $f.$  Hence 
$$ 
\split
f_{11} &= -(n-1)\mu\\
&= -H\,f_1
\endsplit
\tag 4.5
$$ 
where
$H$ is the mean curvature of the level set of $f.$    Note that since $e_1 = \frac{ \n f}{|\n f|}$,  which is a globally defined vector field, we have
$$
f_1 = |\n f|.
$$
The fact that $f$ and $|\n f|$ has the same level sets implies that there exists a function $\beta,$ such that,
$$
f_1 = \beta(f).
$$
In particular,
$$
\split
f_{11} &= e_1e_1 f - \n_{e_1} e_1 f\\
&= \beta'(f)\, f_1 - \la \n_{e_1} e_1, \n f \ra\\
&= \beta'(f)\, \beta(f) - |\n f| \la \n_{e_1} e_1, e_1 \ra\\
& = \beta'(f)\, \beta(f),
\endsplit
$$
hence $f_{11}$ has constant value along the level set of $f$. Combining with \thetag{4.5}, we conclude that the level set of $f$ has constant mean curvature $H$.  
In particular, together with \thetag {4.4}, this implies that the second fundamental form is given by
$$
h_{\alpha \beta} = \frac{H}{(n-1)|\n f|}\, \delta_{\alpha \beta},
$$
a constant multiple of the identity matrix along each level set of $f.$  This implies the splitting of the metric given by the form
$$
ds_M^2 = dt^2 + \eta^2(t)\, ds_N^2,
$$
with
$$
\split
(n-1)\frac{\eta'}{\eta} &= H\\
&= - \frac{f_{11}}{f_1}.
\endsplit
$$
Hence
$$
\eta^{n-1} = C_1\,f_1^{-1},
$$
and
$$
f = C_1\,\int_{0}^t \eta^{-(n-1)} \, dt + C_2
$$
for some constants $C_1$ and $C_2.$

In particular, $g = |\n f|$ implies that
$$
g =C_1\, \eta^{-(n-1)}(t).
$$

The equation 
$$
\d g = -(n-1)\tau\, g + \frac{|\n g|^2}{(n-1) g}
$$
 asserts that
$$
\frac{d^2 g}{d t^2} + (n-1) \eta^{-1} \frac{d \eta}{d t}\, \frac{dg}{d t} = -(n-1)\tau\, g + \frac{|\n g |^2}{(n-1) g},
$$
hence
$$
\frac{d^2 \eta}{d t^2} =  \tau\, \eta.
$$
This implies that $\tau(t)$ must be a function of $t$ alone.\qed
\enddemo

We remark that the manifold $N$ is necessarily compact if $M$ has more than one end.

\proclaim{Corollary 4.2} Let $M^n$ be a complete Riemannian manifold of dimension $n \ge 3.$  Suppose  $M$ satisfies 
the weighted Poincar\'e inequality  
for some nonnegative function $\rho,$ and suppose that the Ricci curvature of $M$ is bounded from below by
$$ 
\text{Ric}_M(x) \ge -\frac{n-1}{n-2}\, \rho(x)
$$ 
for all $x \in M.$ If the volume growth of $M$ satisfies
$$
V_p(R) \le C\, R^{2(n-1)},
$$
then either
\roster
\item $M$ has only one nonparabolic end; or
\item $M$ is given by $M = \Bbb R \times N$ with the warped product metric
$$
ds_M^2 = dt^2 + \eta^2(t)\,ds_N^2,
$$
for some positive function $\eta(t)$, and some compact manifold $N.$ Moreover, $\rho(t)$ is a function of $t$ alone satisfying
$$
(n-2)\eta''\, \eta^{-1} = \rho.
$$
\endroster
\endproclaim

\demo{Proof} According to the theory of Li-Tam \cite{L-T2}, if $M$ has more than one nonparabolic end, then one can construct a 
non-constant bounded harmonic function $f$ with finite Dirichlet integral.  On the other hand, applying Schwarz inequality and using the volume growth assumption, we have
$$
\split
\int_{B(2R)\- B(R)} |\n f|^{\frac{2(n-2)}{n-1}} & \le \left(\int_{B(2R)\-B(R)} |\n f|^2\right)^{\frac{n-2}{n-1}} \left( V_p^{\frac 1{n-1}}(2R) \right)\\
& \le C\, R^2\, \left(\int_{B(2R)\-B(R)} |\n f|^2 \right)^{\frac{n-2}{n-1}}.
\endsplit \tag 4.6
$$
The fact that $f$ has finite Dirichlet integral implies that the right hand side is $o(R^2).$  

Let us denote 
$$
g = |\n f|^{\frac{n-2}{n-1}}.
$$ 
Then Lemma 4.1 and the lower bound of the Ricci curvature assert that 
$$ 
\d g \ge -\rho\, g. \tag 4.7
$$ 
Let $\phi$ be a
non-negative compactly supported smooth function on $M$. Then 
$$ 
\int_M |\n (\phi\,g)|^2= \int_M |\n \phi|^2\,g^2 + 2 \int_M \phi\,g\,\la
\n \phi, \n g\ra +\int_M \phi^2\,|\n g|^2.\tag 4.8
$$ 
The second
term on the right hand side can be written as 
$$ 
\split 
2  \int_M
\phi\,g\,\la \n \phi, \n g\ra &= \frac 12\int_M \la \n
(\phi^2),\n (g^2)\ra\\ & = - \int_M  \phi^2\,g\,\d g - \int_M
\phi^2\,|\n g|^2\\ &=  \int_M \phi^2\,\rho\,g^2- \int_M
\phi^2\,|\n g|^2 -\int_M \phi^2\,g\,(\d g+\rho \,g).
\endsplit
$$
Combining with \thetag {4.8} and property ($\Cal P_\rho$),  
this implies that
$$
\split
\int_M \phi^2\,\rho\,g^2 &\le
 \int_M |\n (\phi\,g)|^2\\
&=\int_M \phi^2\,\rho\,g^2 + \int_M |\n
\phi|^2\,g^2\\
&\qquad -\int_M \phi^2\,g\,\left(\d g+\rho\,g\right).
\endsplit
$$
Hence, we have
$$
\int_M \phi^2\,g\,\left(\d g+\rho \,g\right)\le \int_M |\n
\phi|^2\,g^2. \tag 4.9
$$
For $R_i > 0$, let us  choose $\phi$ to satisfy the properties that
$$
\phi(r(x)) = \left\{ \aligned  1 \qquad &\text{if} \qquad r \le R_i \\
0  \qquad &\text{if} \qquad r \ge 2R_i
\endaligned \right.
$$
and
$$
|\phi'| \le C\,R_i^{-1} \qquad \text{if} \qquad R_i \le r \le 2R_i
$$
for some constant $C>0.$  Then the right hand side of \thetag{4.9}
can be estimated by
$$
\int_M |\n \phi|^2 \, g^2 \le C\,R_i^{-2}\,\int_{B(2R_i) \setminus B(R_i)} g^2.
$$
By the growth estimate of $g$ in \thetag{4.6}, the right hand side tends to 0 as $R_i\to \infty.$  Hence we conclude that \thetag{4.7} must indeed be equality.  
The  corollary now follows from the equality part of Lemma 4.1. Moreover, the manifold $N$ must be compact  because $M$ is assumed to have two ends. \qed
\enddemo

We remark that if the Ricci curvature bound in Corollary 4.2 instead satisfies

$$
\text{Ric}_M(x) \ge -\rho(x),
$$ 
then the validity of the generalized Poincar\'e inequality with weight function $\rho(x)$ alone implies
that every harmonic function of finite Dirichlet energy on $M$ must be constant. In particular, $M$ has
only one nonparabolic end. Indeed, in this case, we have

$$
\d |\n u|\ge -\rho\,|\n u|
$$
for any harmonic function $u.$ With $|\n u|\in L^2(M),$ the inequality must be an equality by a similar argument
using the generalized Poincar\'e inequality. Going back to the Bochner formula, one sees then $|\n u|$ is a constant. 
Since $M$ is nonparabolic, the volume of $M$ must be infinite. The fact that $u$ has finite Dirichlet energy forces 
$|\n u|=0.$ So $u$ is constant. 

This remark is applicable to the case when $M$ is a stable minimal hypersurface in a nonnegatively curved complete manifold, which
recovers a result proved by Schoen and Yau \cite{S-Y1}. Another case is when $M$ is a locally conformally flat manifold with 
scalar curvature $R\le 0.$ Then, according to \cite{S-Y2}, the following inequality holds if $M$ is simply connected.

$$
C(n)\,\left(\int_M |\phi|^{\frac{2n}{n-2}}\right)^{\frac{n-2}{n}}+\int_M \frac{n-2}{4(n-1)}|R| \phi^2\le \int_M |\n \phi|^2
$$
for all $\phi\in C_c^\infty(M),$ where $n=\dim M\ge 3$ and $C(n)>0$ is a constant depending on $n.$
Note that, in particular, this implies that $M$ satisfies a generalized Poincar\'e inequality with weight function 

$$
\rho= \frac{n-2}{4(n-1)}|R|.\tag{4.10}
$$
Also $M$ satisfies a Sobolev inequality, hence all the ends of $M$ must be nonparabolic by \cite{C-S-Z}. In conlcusion, $M$
has only one end in this case provided that the Ricci curvature of $M$ satisfies

$$
\text{Ric}_M(x)\ge -\rho(x)=\frac{n-2}{4(n-1)}R.\tag{4.11}
$$
However, this condition is meaningful only for $n\ge 6.$ Compare this with Corollary 5.5 in the next section.

\proclaim{Corollary 4.3} Let $M^n$ be a complete Riemannian manifold of dimension $n \ge 3.$  Suppose  $M$ satisfies property ($\Cal P_\rho$) 
for some nonnegative function $\rho,$ and suppose that the Ricci curvature of $M$ is bounded from below by
$$ 
\text{Ric}_M(x) \ge -\frac{n-1}{n-2}\, \rho(x)
$$ 
for all $x \in M.$ If $M$ admits a nonconstant harmonic function $f$ with growth estimate satisfying
$$
\int_{B_\rho (2R_i) \- B_\rho(R_i)} \rho\,|\n f|^{\frac{2(n-2)}{n-1}} = o(R_i^2)
$$
for a sequence of $R_i \to \infty,$ then $M$ is given by $M = \Bbb R \times N$ with the warped product metric
$$
ds_M^2 = dt^2 + \eta^2(t)\,ds_N^2,
$$
for some positive function $\eta(t)$, and some manifold $N.$ Moreover, $\rho(t)$ is a function of $t$ alone satisfying
$$
(n-2)\eta''\, \eta^{-1} = \rho.
$$
\endproclaim

\demo{Proof} We use a similar argument as in Corollary 4.2, except that we choose the cut-off function $\phi$ to satisfy the properties that
$$
\phi(r_\rho(x)) = \left\{ \aligned  1 \qquad &\text{if} \qquad r_\rho \le R_i \\
0  \qquad &\text{if} \qquad r_\rho \ge 2R_i
\endaligned \right.
$$
and
$$
|\phi'| \le C\,R_i^{-1} \qquad \text{if} \qquad R_i \le r_\rho \le 2R_i
$$
for some constant $C>0.$ The right hand side of \thetag{4.9} can now be estimated by
$$
\int_M |\n \phi|^2 \, g^2 \le C\,R_i^{-2}\,\int_{B_\rho(2R_i) \setminus B_\rho(R_i)} \rho\, g^2.
$$
The assumption on the growth rate of $g = |\n f|^{\frac{n-2}{n-1}}$ implies that this tends to 0 as $R_i \to \infty.$  
Hence the left hand side of \thetag{4.9} must be identically 0, and Lemma 4.1 implies the corollary.\qed
\enddemo

\heading
\S5 Rigidity and Nonparabolic Ends
\endheading

In this section, we consider a complete manifold, $M^n$, with property ($\Cal P_\rho$) for some nonzero weight function $\rho(x) \ge 0$ 
for all $x \in M.$  Note that since the existence of the weight function $\rho$ is equivalent to $M$ being nonparabolic, $M$ must have 
at least one nonparabolic end.  Assuming that $M$ has at least two nonparabolic ends, $E_1$ and $E_2,$ then a construction of Li-Tam \cite{L-T2} 
asserts that one can construct a nonconstant bounded harmonic function with finite Dirichlet integral.  Indeed, the harmonic function $f$ can be 
constructed by taking a convergent subsequence of the harmonic functions $f_R$, as $R \to \infty$, satisfying
$$
\d f_R = 0 \qquad \text{on} \qquad B(R)
$$
with boundary conditions
$$
f_R = 1 \qquad \text{on} \qquad \p B(R) \cap E_1
$$
and
$$
f_R = 0 \quad \text{on} \qquad \p B(R) \- E_1.
$$
Moreover, the maximum principle asserts that $0 \le f_R \le 1$ for all $R$, hence $0 \le f \le 1.$  We will first prove a lemma concerning the function $f.$

\proclaim{Lemma 5.1} Let $M$ be a complete manifold with property ($\Cal P_\rho$). 
Let $f$ be a bounded harmonic function described above (also in Corollary 2.3) with the property that $0\le \inf f < \sup f < \infty.$  
Let us denote the level set of $f$ at $t$ by
$$
\ell(t) = \{ x \in M \,|\, f(x) =t\}
$$
for $\inf f < t < \sup f$ and we denote the set
$$
\Cal L(a, b) = \{ x \in M\,|\, a< f(x) < b\}
$$
for $ \inf f < a< b< \sup f$.  Then
$$
\int_{\Cal L(a, b)} |\n f|^2 =(b-a) \int_{\ell (b)}|\n f|.
$$
Moreover, 
$$
\int_{\ell(b)} |\n f| = \int_{\ell(t)} |\n f|
$$
for all $\inf f < t < \sup f.$
\endproclaim

\demo{Proof} Let us first observe that if $\phi$ is a nonnegative compactly supported function, then
the co-area formula asserts that
$$
\int_M \phi\, |\n f|^2 \, dV= \int_{\inf f}^{\sup f} \int_{\ell(t)} \phi\, |\n f|\, dA\, dt.
$$
Letting $\phi$ tend to the constant function 1 and using the fact that $f$ has finite Dirichlet integral, we conclude that
$$
\int_{\inf f}^{\sup f} \int_{\ell(t)} |\n f|\, dA\, dt < \infty.
$$
In particular, we conclude that 
$$
\int_{\ell(t)} |\n f|\, dA < \infty
$$
for almost all $\inf f< t < \sup f.$

Let us again denote $\phi$ as a nonnegative compactly supported function. Integrating by parts and using the fact that $f$ is harmonic, we obtain
$$
\split
\int_{\Cal L(a, b)} \phi^2 \, |\n f|^2 &= \int_{\ell(b)} \phi^2\, f\, f_\nu - \int_{\ell(a)} \phi^2\, f\, f_\nu - 2\int_{\Cal L(a, b)}\phi\, f\, \la \n \phi, \n f\ra\\
& = b \int_{\ell(b)} \phi^2\, f_\nu - a\int_{\ell(a)} \phi^2\,  f_\nu -2 \int_{\Cal L(a, b)} \phi\,f\, \la \n \phi, \n f \ra,
\endsplit
$$
where $\nu$ is the unit normal of $\ell(t)$ given by $\nu |\n f| = \n f.$  Since $f_\nu = |\n f|,$ we can write
$$
\int_{\Cal L(a, b)} \phi^2\, |\n f|^2 = b \int_{\ell(b)} \phi^2\, |\n f| - a\int_{\ell(a)} \phi^2\, |\n f| -2 \int_{\Cal L(a, b)} \phi\,f\, \la \n \phi, \n f \ra.\tag 5.1
$$
Let us choose
$$
\phi(x) = \left\{ \aligned 1 \quad \qquad & \qquad \text{on} \qquad B_\rho(R)\\
R+1- r_\rho(x) & \qquad \text{on} \qquad B_\rho(R+1) \- B_\rho(R)\\
0\quad \qquad & \qquad \text{on} \qquad M \- B_\rho(R+1). \endaligned \right.
$$
Then
$$
|\n \phi| = \left\{ \aligned \sqrt{\rho} & \qquad \text{on} \qquad B_\rho(R+1)\- B_\rho(R)\\
0 \quad & \qquad \text{otherwise} \endaligned \right.
$$
The last term on the right hand side of \thetag{5.1} can be estimated by
$$
\split
&\left|\int_{\Cal L(a, b)}\phi\,f\, \la \n \phi, \n f \ra \right| \\
&\quad \qquad \le \int_{\Cal L(a, b) \cap (B_\rho(R+1) \- B_\rho(R))} \sqrt{\rho}\,|f|\,  |\n f|\\
&\quad \qquad \le \left(\int_{\Cal L(a, b) \cap (B_\rho(R+1) \- B_\rho(R))} \rho\,f^2 \right)^{\frac 12} 
\left( \int_{\Cal L(a, b) \cap (B_\rho(R+1) \- B_\rho(R))} |\n f|^2\right)^{\frac12}\\
&\quad \qquad \le C\, \exp(-R)\,\left(\int_{\Cal L(a, b) \cap (B_\rho(R+1) \- B_\rho(R))} \rho\,f^2 \right)^{\frac 12},
\endsplit \tag5.2
$$
where we have used the estimate provided by Corollary 2.3.

Following the notation and the estimate of Corollary 2.3, let $E_1$ be the nonparabolic end on which $\sup f $ is achieved at infinity.  
Then on the set $(B_\rho(R+1) \- B_\rho(R))\- E_1$, inequality \thetag{2.9} implies that
$$
\int_{(B_\rho(R+1) \- B_\rho(R))\-E_1} \rho\,(f - \inf f)^2 \le C\, \exp(-2R).
$$
In particular, 
$$
\split
&\int_{(\Cal L(a, b) \cap (B_\rho(R+1) \- B_\rho(R))) \- E_1} \rho\,f^2 \\
& \qquad \le b^2\, \int_{(\Cal L(a, b) \cap (B_\rho(R+1) \- B_\rho(R))) \- E_1} \rho\\
& \qquad \le b^2\, (a- \inf f)^{-2}\,\int_{(\Cal L(a, b) \cap (B_\rho(R+1) \- B_\rho(R))) \- E_1} \rho\,(f - \inf f)^2\\
&\qquad \le b^2\, (a - \inf f)^{-2}\, C\, \exp (-2R).
\endsplit \tag5.3
$$
Similarly, on $E_1$, \thetag {2.10} implies that we have
$$
\split
&\int_{\Cal L(a, b) \cap (B_\rho(R+1) \- B_\rho(R))\cap E_1} \rho\,f^2\\
&\qquad \le b^2 \int_{\Cal L(a, b) \cap (B_\rho(R+1) \- B_\rho(R))\cap E_1} \rho\\
&\qquad \le b^2 \,(\sup f - b)^{-2}\int_{(B_\rho(R+1) \- B_\rho(R)) \cap E_1} \rho\, (\sup f-f)^2 \\
&\qquad \le b^2\, (\sup f-b)^{-2}\,C\,\exp(-2R)
\endsplit
$$
Together with \thetag{5.3} and \thetag{5.2}, we conclude that
$$
\left|\int_{\Cal L(a, b)} \phi\,f\, \la \n \phi, \n f \ra \right| \le (b((a-\inf f)^{-1} + (\sup f - b)^{-1}))\, C\,\exp(-2R).
$$
Letting $R \to \infty,$ \thetag{5.1} becomes
$$
\int_{\Cal L(a, b)} |\n f|^2 = b \int_{\ell(b)} |\n f| - a\int_{\ell(a)}  |\n f|.\tag5.4
$$

We now observe that since \thetag {5.4} is independent of $\inf f$ and $\sup f$,  if we apply \thetag{5.4} to the function  
$f+\epsilon,$  then we have
$$
\int_{\Cal L(a, b)} |\n f|^2 = (b+\epsilon) \int_{\ell(b)} |\n f| - (a +\epsilon)\int_{\ell(a)}  |\n f|.
$$
Combining with \thetag{5.4}, we conclude that
$$
\int_{\ell(a)} |\n f| = \int_{\ell(b)} |\n f|.
$$
Since $a$ is arbitrary, this proves the lemma.\qed
\enddemo

We are now ready to prove the first main theorem. Let us first define
$$
S(R) = \sup_{B_\rho(R)} \sqrt{\rho}
$$
to be the supremum of $\sqrt{\rho}$ over the set $B_\rho(R).$

\proclaim{Theorem 5.2} Let $M^n$ be a complete manifold with dimension $n \ge 3.$  Assume that $M$ satisfies property ($\Cal P_\rho$) 
for some nonzero weight function $\rho \ge 0.$  Suppose the Ricci curvature of $M$ satisfies the lower bound
$$
\text{Ric}_M(x) \ge - \frac{n-1}{n-2}\, \rho(x)
$$
for all $x \in M.$  If $\rho$ satisfies the growth estimate

$$
\liminf_{R\to \infty} \frac{S(R)}{F(R)}=0,
$$
where

$$
F(R) = \left\{ \aligned \exp \left(\frac{n-3}{n-2}\,R\right) & \qquad \text{when} \qquad n\ge 4\\
R \quad \qquad & \qquad \text{when} \qquad n=3,\endaligned \right.
$$
then either
\roster
\item $M$ has only one nonparabolic end; or
\item $M$ has two nonparabolic ends and is given by $M = \Bbb R \times N$ with the warped product metric
$$
ds_M^2 = dt^2 + \eta^2(t)\,ds_N^2,
$$
for some positive function $\eta(t)$, and some compact manifold $N.$ Moreover, $\rho(t)$ is a function of $t$ alone satisfying
$$
(n-2)\eta''\, \eta^{-1} = \rho.
$$
\endroster
\endproclaim

\demo{Proof} Let us assume that $M$ has at least two nonparabolic ends. Then there exists a bounded harmonic 
function $f$ with finite Dirichlet integral constructed as above. We may assume that $\inf f = 0$ and $\sup f = 1.$
Note that the improved Bochner formula asserts that
$$
\d g \ge -\rho \,g
$$
where $g = |\n f|^{\frac{n-2}{n-1}},$ and according to Lemma 4.1, it suffices to show that
$$
\d g = -\rho \, g.
$$

To see this, let us consider $\phi$ to be a non-negative  smooth function with compact support in $M$. Then 
$$ 
\int_M |\n (\phi\,g)|^2= \int_M |\n \phi|^2\,g^2 + 2 \int_M \phi\,g\,\la
\n \phi, \n g\ra +\int_M \phi^2\,|\n g|^2.\tag 5.5
$$ 
The second term on the right hand side can be written as 
$$ 
\split 
2  \int_M
\phi\,g\,\la \n \phi, \n g\ra &= \frac 12\int_M \la \n
(\phi^2),\n (g^2)\ra\\ & = - \int_M  \phi^2\,g\,\d g - \int_M
\phi^2\,|\n g|^2\\ &=  \int_M \phi^2\,\rho\,g^2- \int_M
\phi^2\,|\n g|^2 -\int_M \phi^2\,g\,(\d g+\rho \,g).
\endsplit
$$
Combining with \thetag {5.5} and property ($\Cal P_\rho$),  
this implies that
$$
\split
\int_M \phi^2\,\rho\,g^2 &\le
 \int_M |\n (\phi\,g)|^2\\
&=\int_M \phi^2\,\rho\,g^2 + \int_M |\n
\phi|^2\,g^2\\
&\qquad -\int_M \phi^2\,g\,\left(\d g+\rho\,g\right).
\endsplit
$$
Hence, we have
$$
\int_M \phi^2\,g\,\left(\d g+\rho \,g\right)\le \int_M |\n
\phi|^2\,g^2. \tag 5.6
$$

Let us choose $\phi = \psi\,\chi$ to be the product of two compactly supported functions.  For $0< \delta <1$ and $0< \epsilon < \frac 12$, let us choose $ \chi$ to be
$$
\chi(x) = \left\{ \aligned 0 \qquad \qquad \quad & \qquad \text{on} \qquad \Cal L(0, \delta \epsilon) \cup \Cal L(1-\delta\epsilon , 1) \\
(-\log \delta)^{-1} (\log f - \log (\delta \epsilon))  & \qquad \text{on} \qquad \Cal L(\delta\epsilon, \epsilon) \cap (M\- E_1)\\
(-\log \delta)^{-1} (\log (1-f) - \log (\delta \epsilon)) & \qquad \text{on} \qquad \Cal L(1-\epsilon, 1- \delta\epsilon) \cap E_1\\
1 \qquad \qquad \quad &  \qquad \text{otherwise.} 
\endaligned \right.
$$
For $R >0$, we choose
$$
\psi(x) = \left\{ \aligned 1\quad & \qquad \text{on} \qquad B_\rho(R-1)\\
R-r_\rho  & \qquad \text{on} \qquad B_\rho(R) \- B_\rho(R-1)\\
0 \quad & \qquad \text{on} \qquad M\-B_\rho(R) \endaligned \right.
$$
Then applying to the right hand side of \thetag{5.6},
we obtain
$$
\int_M |\n \phi|^2 \, g^2 \le 2\int_M |\n \psi|^2\, \chi^2 \, |\n f|^{\frac{2(n-2)}{n-1}} 
+ 2\int_M |\n \chi|^2 \, \psi^2 \, |\n f|^{\frac{2(n-2)}{n-1}}.\tag 5.7
$$

Let us now recall that, using the assumption on the Ricci curvature,  the  local gradient estimate of Cheng-Yau \cite{C-Y} (see \cite{L-W2}) 
for positive harmonic functions asserts that for all $R_0>0$, 
$$
|\n f|(x) \le \left((n-1)\sup_{B(x, R_0)}\sqrt{\rho(y)} + C\,R_0^{-1} \right) \,f(x),\tag 5.8
$$
where $C$ is a constant depending only on $n,$ and $B(x, R_0)$ is the ball of radius $R_0$ centered at $x$ with respect to the background metric $ds_M^2.$ 
Let us now choose $R_0 = (\sup_{B(x, R_0)} \sqrt{\rho})^{-1}.$  This choice of $R_0$ is possible as the function $r - (\sup_{B(x, r)} \sqrt{\rho})^{-1}$ 
is negative when $r \to 0$ and it tends to $\infty$ as $r \to \infty.$
Let us observe that if $y \in B(x, R_0)$,  and if $\gamma$ is a $ds_M^2$ minimizing geodesic joining $x$ to $y$, then
$$
\split
r_\rho(x, y) &= \int_{\gamma} \sqrt{\rho(\gamma(t)}\, dt\\
&\le \sup_{B(x, R_0)} \sqrt{\rho(y)}\, R_0\\
&\le 1.
\endsplit
$$
This implies that $B(x, R_0) \subset B_\rho(x,1).$ Hence \thetag{5.8} can be written as
$$
|\n f|(x) \le C\, (\sup_{B_\rho(x, 1)} \sqrt{\rho})\, f(x). \tag 5.9
$$
Similarly, applying the same estimate to $1-f$, we also have
$$
|\n f|(x) \le C\, (\sup_{B_\rho(x, 1)} \sqrt{\rho})\, (1-f(x)).\tag 5.10
$$

At the end $E_1$, the first term on the right hand side of \thetag{5.7} can be estimated by
$$
\split
\int_{E_1} |\n \psi|^2\, \chi^2 \, |\n f|^{\frac{2(n-2}{n-1}} &\le \int_{\Omega} \rho\,|\n f|^{\frac {2(n-2)}{n-1}}\\
&\le \left(\int_{\Omega} |\n f|^2 \right)^{\frac{n-2}{n-1}} \left( \int_{\Omega} \rho^{n-1}\right)^{\frac 1{n-1}},
\endsplit \tag5.11
$$
where
$$
\Omega = E_1 \cap (B_\rho(R) \- B_\rho(R-1)) \cap \Cal L(\delta\epsilon, 1-\delta\epsilon).
$$
Applying Corollary 2.3, we conclude that
$$
\left(\int_{\Omega} |\n f|^2 \right)^{\frac{n-2}{n-1}} \le C\, \exp\left(-\frac{2(n-2)R}{n-1}\right).\tag5.12
$$
On the other hand, using \thetag{2.10}, we have
$$
\split
\int_{\Omega} \rho^{n-1}
& \le S^{2(n-2)}(R) \,  \int_{\Omega} \rho\\
 &   \le S^{2(n-2)}(R) \, (\delta\epsilon)^{-2}\, \int_{\Omega} \rho\, (1-f)^2\\
 & \le C\, S^{2(n-2)}(R) \,  (\delta\epsilon)^{-2}\, \exp(-2R).
 \endsplit
 $$
 Hence together with \thetag{5.11} and \thetag{5.12}, we obtain
 $$
 \int_{E_1} |\n \psi|^2\, \chi^2 \, |\n f|^{\frac{2(n-2)}{n-1}} \le C\,S^{\frac{2(n-2)}{n-1}}(R) \, (\delta\epsilon)^{-\frac 2{n-1}}\, \exp(-2R).\tag 5.13
 $$

Using \thetag{5.10}, the second term on the right hand side of \thetag{5.7} at $E_1$ can be estimated by
$$
\split
\int_{E_1}& |\n \chi|^2 \, \psi^2 \, |\n f|^{\frac{2(n-2)}{n-1}}\\ & \le (\log \delta)^{-2}\,
\int_{\Cal L(1-\epsilon, 1-\delta\epsilon) \cap E_1 \cap B_\rho(R)} |\n f|^{2+ \frac{2(n-2)}{n-1}}\,(1-f)^{-2}\\
&\le C\, S^{\frac{2(n-2)}{n-1}}(R+1)\, (\log \delta)^{-2}\,\int_{\Cal L(1-\epsilon, 1-\delta\epsilon) \cap E_1 \cap B_\rho(R)} |\n f|^2\,(1-f)^{\frac{2(n-2)}{n-1}-2}.
\endsplit\tag5.14
$$
Note that the co-area formula and Lemma 5.1 imply that
$$
\split
\int_{\Cal L(1-\epsilon, 1-\delta\epsilon) \cap E_1 \cap B_\rho(R)}& |\n f|^2\,(1-f)^{\frac{2(n-2)}{n-1}-2}\\
&\le \int_{1-\epsilon}^{1-\delta \epsilon} (1-t)^{\frac{2(n-2)}{n-1}-2} \int_{\ell(t) \cap E_1 \cap B_\rho(R)} |\n f|\,dA\,dt\\
&\le \int_{\ell(b)} |\n f|\,dA \int_{1-\epsilon}^{1-\delta \epsilon} (1-t)^{\frac{2(n-2)}{n-1} -2}\, dt
\endsplit
$$
for any level $b$.
Since 
$$
\int_{1-\epsilon}^{1-\delta \epsilon} (1-t)^{\frac{2(n-2)}{n-1} -2}\, dt 
 = \left\{ \aligned \frac{n-1}{n-3}\,\left(\epsilon^{\frac{n-3}{n-1}} - (\delta \epsilon)^{\frac{n-3}{n-1}}\right)& \qquad \text{when} \qquad n\ge 4\\
- \log \delta \qquad \qquad &\qquad \text{when} \qquad n=3,
\endaligned \right.
$$
together with \thetag{5.14}, we conclude that
$$
\split
\int_{E_1} &|\n \chi|^2 \, \psi^2 \, |\n f|^{\frac{2(n-2)}{n-1}} \\
& \le \left\{ \aligned C\, S^{\frac{2(n-2)}{n-1}}(R+1)\,(\log \delta)^{-2}\, \left( 1- \delta^{\frac{n-3}{n-1}}\right)\,\epsilon^{\frac{n-3}{n-1}}
& \qquad \text{when} \qquad n\ge 4\\
C\, S(R+1)\,(-\log \delta)^{-1} \qquad \qquad & \qquad \text{when} \qquad n=3.
\endaligned \right.
\endsplit \tag5.15
$$

When $n\ge 4,$ we set $\delta = \frac 12$, \thetag{5.7}, \thetag{5.13}, and \thetag{5.15} together imply that
$$
\int_{E_1} |\n \phi|^2 \, g^2 \le C\,S^{\frac{2(n-2)}{n-1}}(R+1)\,\left( \epsilon^{-\frac 2{n-1}}\, \exp(-2R) + \epsilon^{\frac{n-3}{n-1}}\right).
$$
A similar argument using the function $f$ instead of $1-f$ and \thetag{5.9} instead of \thetag{5.10} will yield the estimate
$$
\int_{M\-E_1} |\n \phi|^2 \, g^2 \le C\,S^{\frac{2(n-2)}{n-1}}(R+1)\,\left( \epsilon^{-\frac 2{n-1}}\, \exp(-2R) + \epsilon^{\frac{n-3}{n-1}}\right)
$$
on the other ends of $M.$
In particular, we conclude that
$$
\int_M |\n \phi|^2 \, g^2 \le C\,S^{\frac{2(n-2)}{n-1}}(R+1)\, \left(\epsilon^{-\frac 2{n-1}}\, \exp(-2R) + \epsilon^{\frac{n-3}{n-1}}\right)
$$
is valid for $R$ sufficienly large.  Setting $\epsilon = \exp(-2R),$ this becomes
$$
\int_M |\n \phi|^2 \, g^2 \le C\,S^{\frac{2(n-2)}{n-1}}(R+1)\,  \exp\left(-\frac{2(n-3)R}{n-1} \right).
$$
Hence the right hand side tends to 0 by taking $R+1= R_i$ 
and letting $i \to \infty,$ where 

$$
\lim_{R_i\to \infty}\frac{S(R_i)}{F(R_i)}=0,
$$
and the theorem follows for $n \ge 4.$

Similarly, when $n=3$, we get

$$
\int_M |\n \phi|^2 \, g^2 \le C\,S(R+1)\,\left( \delta^{-1} \epsilon^{-1}\, \exp(-2R) + (-\log \delta)^{-1}\right).\tag5.16
$$
Set $\delta = \epsilon = \exp(-R\,q(R))$ with
$$
q(R) = \sqrt{\frac{S(R+1)}{R}}.
$$
Note that the assumption that $\liminf \frac{S(R)}{R} = 0$ implies that $q(R_i) \to 0$ for some sequence of $R_i \to \infty.$  Hence
$$
\split
S(R_i+1)\, & \left( \delta^{-1} \epsilon^{-1} \exp(-2R_i) + (-\log \delta)^{-1} \right) \\
&= S(R_i+1) \left( \exp(-2R_i(1-q(R_i))) + R_i^{-1} q^{-1}(R_i) \right).
\endsplit
$$
For sufficiently large $i$, the term $q(R_i) \le \frac 12,$ and 
$$
S(R_i+1)\, \exp(-2R_i(1- q(R_i))) \le S(R_i+1) \, \exp(-R_i)
$$\
must tend to $0$ as $i \to \infty.$  By definition of $q$, the other term 
$$
S(R_i+1)\, R_i^{-1} q^{-1}(R_i) = q(R_i)
$$
also tends to 0 as remarked earlier.  Hence  the right hand side of \thetag{5.16} tends to $0$  and the theorem is proved.\qed
\enddemo

Obviously, if $\rho$ is bounded, then the growth assumption on $\rho$ in the theorem is automatically satisfied. We claim this is also true
when $\rho(r)$ is an non-decreasing function depending only on the distance to a compact set with respect to the background metric $ds_M^2$ and $n\ge 4.$
Indeed, the $\rho$-distance to the compact set

$$
r_\rho(r) = \int_0^r \sqrt{\rho(t)}\, dt
$$
is also a function of $r$ alone and the function 

$$
h(r)=\exp \left(-\frac{n-3}{n-2}\,r_\rho(r)\right)\to 0
$$ 
as $r\to \infty.$ So there exists $r_i\to \infty$ such that

$$
\lim_{r_i\to \infty} h'(r_i)=0,
$$ 
which verifies the claim. 

When $n=3,$ if we assume $\rho(r)$ is non-decreasing and of subexponential growth in $r,$
then the same argument applying to the function $\ln r_\rho(r)$ instead will show the growth assumption of the theorem is satisfied.

\proclaim{Corollary 5.3} Let $M^n$ be a complete manifold, of dimension $n \ge 4$, with property ($\Cal P_\rho$) for some nonzero weight function $\rho \ge 0.$  
Suppose the Ricci curvature of $M$ satisfies the lower bound
$$
\text{Ric}_M(x) \ge - \frac{n-1}{n-2}\, \rho(x)
$$
for all $x \in M.$  Assume that $\rho(x)=\rho(r(x))$ is a non-decreasing function of the distance $r$ to a fixed compact set with respect to the background metric $ds_M^2.$ 
Then either
\roster
\item $M$ has only one nonparabolic end; or
\item $M$ is given by $M = \Bbb R \times N$ with the warped product metric
$$
ds_M^2 = dt^2 + \eta^2(t)\,ds_N^2,
$$
for some positive function $\eta(t)$, and some compact manifold $N.$ Moreover, $\rho(t)$ is a function of $t$ alone satisfying
$$
(n-2)\eta''\, \eta^{-1} = \rho.
$$
\endroster
In particular, if $\rho(r)$ is asymptotically a polynomial in $r$, then the conclusion of the theorem is valid.
\endproclaim

\proclaim{Corollary 5.4} Let $M^3$ be a complete manifold, of dimension $3$, with property ($\Cal P_\rho$) for some nonzero weight function $\rho \ge 0.$  
Suppose the Ricci curvature of $M$ satisfies the lower bound
$$
\text{Ric}_M(x) \ge - 2 \rho(x)
$$
for all $x \in M.$  Assume that $\rho(x)=\rho(r(x))$ is a non-decreasing function of the distance to a fixed compact set with respect to the background metric $ds_M^2$ 
and that $\rho(r)$ is of subexponential growth. Then either
\roster
\item $M$ has only one nonparabolic end; or
\item $M$ has two nonparabolic ends and is given by $M = \Bbb R \times N$ with the warped product metric
$$
ds_M^2 = dt^2 + \eta^2(t)\,ds_N^2,
$$
for some positive function $\eta(t)$, and some compact manifold $N.$ Moreover, $\rho(t)$ is a function of $t$ alone satisfying
$$
\eta''\, \eta^{-1} = \rho.
$$
\endroster
In particular, if $\rho(r)$ is asymptotically a polynomial in $r$, then the conclusion of the theorem is valid.
\endproclaim

When $M$ is an $n$-dimensional, locally comformally flat, simply connected manifold with nonpositive scalar curvature $R$, according to the remarks
made after Corollary 4.2, all ends of $M$ are nonparabolic and $M$ satisfies a generalized Poincar\'e inequality
with weight function $\rho$ given by \thetag{4.10}. Applying Theorem 5.2, we have the following result.

\proclaim{Corollary 5.5} Let $(M^n, ds^2_M)$ be a complete, simply connected, locally comformally flat manifold, of dimension $n \ge 4$, with nonpositive scalar curvature $R.$
Suppose the Ricci curvature of $M$ satisfies the lower bound
$$
\text{Ric}_M(x) \ge  \frac{1}{4}\, R(x)
$$
for all $x \in M.$  Assume that the conformal metric $|R|\,ds^2_M$ is also complete, and 
$$
\liminf_{r\to \infty} \exp\left(-\frac{n-3}{n-2}r\right)\,\sup_{B_{|R|}(r)} |R|(x)=0,
$$ 
where the supremum is taken over the ball of radius $r$, with respect to the metric $|R|\,ds^2_M$, centered at some fixed point $p \in M.$
Then either
\roster
\item $M$ has only one end; or
\item $M$ is given by $M = \Bbb R \times N$ with the warped product metric
$$
ds_M^2 = dt^2 + \eta^2(t)\,ds_N^2,
$$
for some positive function $\eta(t)$, and some compact manifold $N.$ Moreover, $R$ is a function of $t$ alone satisfying
$$
4(n-1)\eta''\, \eta^{-1} =-R.
$$
\endroster
\endproclaim

\heading
\S6 Warped Product Metrics 
\endheading

Let $M^n = \br \times N^{n-1}$ be the product manifold endowed with the warped product metric
$$
ds_M^2 = dt^2 + \eta^2(t)\, ds_N^2,
$$
where $ds_N^2$ is a given metric on the compact manifold $N.$  Our purpose is to compute the curvature on $M$  
and discuss the possibility of the existence of a weight function $\rho$ that is given by a multiple of the lower bound of the Ricci curvature.

Let $\{\bar \omega_2, \dots,  \bar \omega_n\}$ be an orthonormal coframe on $N$ with respect to $ds_N^2$.  If we define $\omega_1 = dt$ 
and $\omega_\alpha = \eta(t)\, \bar \omega_\alpha$ for $2 \le \alpha \le n,$ then the set $\{\omega_i\}_{i=1}^n$ forms an orthonormal coframe 
of $M$ with respect to $ds_M^2.$  The first structural equations assert that
$$
d \omega_i = \omega_{ij} \wedge \omega_j,
$$
where $\omega_{ij}$ are the connection 1-forms with the property that
$$
\omega_{ij} = -\omega_{ji}.
$$
On the other hand, direct exterior differentiation yields
$$
d \omega_1 = 0
$$
and
$$
\split
 d \omega_\alpha &= \eta' \,\omega_1 \wedge \bar \omega_\alpha  + \eta\,\bar \omega_{\alpha \beta}\, \wedge \bar\omega_{\beta}\\
 &= -(\log \eta)' \,\omega_\alpha \wedge \omega_1 + \bar \omega_{\alpha \beta} \wedge \omega_{\beta},
 \endsplit
 $$
 where $\bar \omega_{\alpha \beta}$ are the connection 1-forms on $N$ and $\eta'$ is the derivative of $\eta$ with respect to $t.$
 Hence we conclude that the connection 1-forms are given by
 $$
 \split
 \omega_{1 \alpha} &= - \omega_{\alpha 1}\\
 & = (\log \eta)' \,\omega_\alpha
 \endsplit \tag 6.1
 $$
 and
 $$ 
 \omega_{\alpha \beta} = \bar \omega_{\alpha \beta}. \tag 6.2
 $$

 The second structural equations also assert that
 $$
 d\omega_{ij} - \omega_{ik} \wedge \omega_{kj} = \frac 12 R_{ijkl}\, \omega_l \wedge \omega_k,
 $$
 where $R_{ijkl}$ is the curvature tensor on $M$.
 Exterior differentiating \thetag{6.1} yields
  $$
 d \omega_{1 \alpha} = (\log \eta)'' \,\omega_1 \wedge \omega_{\alpha} + (\log \eta)' \left(-(\log \eta)' \,\omega_\alpha \wedge \omega_1 + 
\bar \omega_{\alpha \beta} \wedge \omega_{\beta}\right).
 $$
 Hence combining with \thetag{6.1} and \thetag{6.2}, we have
 $$
 d \omega_{1 \alpha} - \omega_{1 \beta} \wedge \omega_{\beta \alpha} =  \left((\log \eta)'' + ((\log \eta)')^2 \right)\,\omega_1 \wedge \omega_{\alpha}.
 $$
 Also, exterior differentiating \thetag{6.2} gives
 $$
 d \omega_{\alpha \beta} = d \bar \omega_{\alpha \beta},
 $$
 and
 $$
 \split
 d \omega_{\alpha \beta} - \omega_{\alpha 1} \wedge \omega_{1 \beta} - \omega_{\alpha \gamma} \wedge \omega_{\gamma \beta} 
&= d \bar \omega_{\alpha \beta} - \bar \omega_{\alpha \gamma} \wedge \bar \omega_{\gamma \beta} + ((\log \eta)')^2 \,\omega_\alpha \wedge \omega_\beta\\
 &= \frac 12 \bar R_{\alpha \beta \gamma \tau} \,\bar \omega_\tau \wedge \bar \omega_\gamma + ((\log \eta)')^2\, \omega_\alpha \wedge \omega_\beta\\
 & =  \frac 12 \bar R_{\alpha \beta \gamma \tau}\,  \eta^{-2}\, \omega_\tau \wedge \omega_\gamma + ((\log \eta)')^2 \,\omega_\alpha \wedge \omega_\beta,
 \endsplit
 $$
 where $\bar R_{\alpha \beta \gamma \tau}$ is the curvature tensor on $N$.  In particular, the sectional curvature of the two-plane section spanned 
by $e_1$ and $e_\alpha$ is given by
$$
K(e_1, e_\alpha) = -\left((\log \eta)'' + ((\log \eta)')^2\right).
$$
Also the sectional curvature of the two plane section spanned by $e_\alpha$ and $e_\beta$ is given by
$$
K(e_\alpha, e_\beta) = \eta^{-2} \,\bar K(e_\alpha, e_\beta) - ((\log \eta)')^2,
$$
where $\bar K$ is the sectional curvature of $N.$  Moreover the curvature tensor is given by
$$
R_{1 \alpha jk} = \left\{ \aligned (\log \eta)'' &+ ((\log \eta)')^2 \qquad \text{if } j=\alpha, \, k=1\\
- (\log \eta)'' &- ((\log \eta)')^2 \qquad \text{if } j=1, \, k=\alpha\\
&0 \qquad \text{otherwise}.\endaligned \right.
$$
and
$$
R_{\alpha \beta ij} =\left\{ \aligned \eta^{-2}\,\bar R_{\alpha \beta \gamma \tau} +& ((\log \eta)')^2\,(\delta_{\alpha \tau} \delta_{\beta \gamma} - 
\delta_{\alpha \gamma}\delta_{\beta \tau}) \qquad \text{if } i=\gamma, j = \tau\\
0& \qquad \text{otherwise}.\endaligned \right.
$$
The Ricci curvature is then given by
$$
\split
\text{Ric}_{1 j} &= \sum_{\alpha} R_{1 \alpha j \alpha} \\
&= -(n-1)\left( (\log \eta)'' + (( \log \eta)')^2\right) \delta_{1 j}\\
&= - (n-1) \eta''\, \eta^{-1}\,\delta_{1j},
\endsplit \tag 6.3
$$
and
$$
\split
\text{Ric}_{\alpha \beta} & = \sum_{\gamma \neq \alpha} R_{\alpha \gamma \beta \gamma} + R_{\alpha 1 \beta 1}\\
&=\eta^{-2}\, \bar{\text{Ric}}_{\alpha \beta} -\left( (\log \eta)'' + (n-1)(( \log \eta)')^2 \right)\,\delta_{\alpha \beta},
\endsplit \tag 6.4
$$
where $\bar{\text{Ric}}_{\alpha \beta}$ is the Ricci tensor on $N.$

Let $\d$ be the Laplace operator defined on functions of $M.$  When restricted to a function
$f(t)$ depending only on the variable $t$, it has the expression
$$
\d f(t) = \frac{\partial^2 f}{\partial t^2} + (n-1) (\log \eta)'\, \frac{\partial f}{\partial t}.
$$
In particular, if we define
$$
f(t) = \int_0^t \eta^{-(n-1)}\,ds, \tag 6.5
$$
then a direct computation yields that
$$
\d f = 0.
$$
Taking 
$$
\split
g &= |\n f|^{\frac{n-2}{n-1}}\\
& = \eta^{-(n-2)},
\endsplit\tag 6.6
$$
we have
$$
\split
\d g &= -(n-2) \eta''\,\eta^{-(n-1)} \\
&= - (n-2)  \eta''\,\eta^{-1}\,g
\endsplit
$$
According to Proposition 1.1, if 
$$
\eta'' > 0,\tag 6.7
$$
then the function
$$
\rho = (n-2) \eta''\, \eta^{-1}
$$
is a weight function, and $M$ satisfies the corresponding weighted Poincar\'e inequality.

Obviously, \thetag{6.3} implies that
$$
\text{Ric}_{1j} = -\frac{n-1}{n-2} \, \rho\, \delta_{1j}.
$$
The condition that
$$
\text{Ric}_M \ge - \frac{n-1}{n-2} \, \rho
$$
is then equivalent to
$$
\eta^{-2} \,\bar{\text{Ric}}_{\alpha \beta}  - ((\log \eta)'' + (n-1)((\log \eta)')^2)\, \delta_{\alpha \beta} \ge - (n-1) \eta^{-1} \, \eta''\, \delta_{\alpha \beta}
$$
because of \thetag {6.4}.  This can be rewritten as
$$
(n-2) (\log \eta)''  +\eta^{-2}\, \text{Ric}_N \ge 0.\tag 6.8
$$

Let us summarize the above computation in the following proposition.

\proclaim{Proposition 6.1} Let us consider the warped-product metric defined on $M^n = \Bbb R \times N^{n-1}$, given by the form
$$
ds^2_M = dt^2 + \eta^2(t)\,ds^2_N
$$
for some positive function $\eta(t)$ depending on $t$ alone.
If the warping function $\eta$ satisfies \thetag {6.7} and \thetag {6.8}, then the warped product manifold $M$ will have the property that
$$
\text{Ric}_M \ge - \frac{n-1}{n-2} \rho,
$$
where
$$
\rho = (n-2) \eta''\,\eta^{-1}.\tag 6.9
$$
Moreover, $M$ will satisfy the weighted Poincar\'e inequality
$$
\int_M \rho\, \phi^2 \le \int_M |\n \phi|^2
$$
for any compactly supported function $\phi \in C_c^{\infty}(M).$  Also
$f$ defined by \thetag {6.5} is a harmonic function with $g$ defined by \thetag {6.6} satisfying
$$
\d g = - \rho\, g,
$$
and
$$
\text{Ric}(\n f, \n f) = - \frac{n-1}{n-2} \rho \, |\n f|^2.
$$
\endproclaim

\noindent
{\bf Example 6.2:} We will provide the following special examples for $\eta$.  Let us consider
$$
\eta(t) = \cosh u(t).
$$
We will take $u(t)$ to be an odd function with the property that
$$
u(-t) = -u(t).
$$
Let us first focus our computation on the positive half-line $[0, \infty).$  Direct differentiation gives
$$
\eta' = u' \,\sinh u
$$
and
$$
\eta'' = u''\, \sinh u + (u')^2 \cosh u.
$$
If we assume that $\text{Ric}_N \ge -C$ for some constant $C \ge 0$,  then the conditions
\thetag {6.7} and \thetag{6.8} becomes
$$
u''\, \sinh u + (u')^2 \cosh u >0
$$
and
$$
\split
(n-2)(\eta'' \eta - (\eta')^2) &= (n-2)(u''\, \sinh u\, \cosh u + (u')^2)\\
&\ge C.
\endsplit
$$
Obviously, if we take any $u$ with the property that 
$$
u'' \ge 0
$$ 
and 
$$
u' > \sqrt{\frac{C}{n-2}}
$$
on $[0, \infty)$,
then both of these conditions will be satisfied.  The oddness of $u$ ensures that the conditions are valid also on $(-\infty, 0].$  
Examples of such $u$ can be given by a smooth approximation of the function
$$
v(t) = \left\{\aligned &C_1\, t^{\alpha} \qquad \text{for } \qquad 1 \le t\\
&C_1\, t \qquad \text{for } \qquad 0 \le t\le 1\endaligned \right.
$$
for any value of $\alpha \ge 1$ and $C_1 > \sqrt{\frac{C}{n-2}}.$  In this case,
$$
\rho(t) \sim (n-2)\, \alpha^2\, C_1^2\, t^{2\alpha -2}
$$
as $|t| \to \infty,$ hence $M$ has property ($\Cal P_\rho$).

On the other hand, we will also show that for some cases of $\rho$, the warped product situation does not exist.

\proclaim{Theorem 6.3}Let $M^n = \Bbb R \times N^{n-1}$ be the warped-product manifold with metric given by
$$
ds^2_M = dt^2 + \eta^2(t)\,ds^2_N.
$$
Suppose $M$ has property ($\Cal P_\rho$)
with
$$
\rho = (n-2)\eta''\,\eta^{-1}
$$
and
$$
\text{Ric}_M \ge - \frac{n-1}{n-2} \rho.
$$
Then on an nonparabolic end $E$, 
$$
\liminf_{x\to \infty} \rho(x) > 0,
$$
for $x \in E.$
\endproclaim

\demo{Proof} Combining \thetag {6.9} and \thetag {6.8}, we have
$$
\frac{\rho}{n-2} \eta^2 - (\eta')^2 \ge -C\tag 6.10
$$
if we assume that $(n-2)^{-1}\text{Ric}_N \ge C$ for some constant $C.$  Let $E$ be a nonparabolic end. 
By reparametrizing if necessary, we may assume that $E$ is given by $[0, \infty) \times N.$ The geometric condition for nonparabolicity is then given by
$$
\int_0^\infty \eta^{-(n-1)}\, dt < \infty.
$$
In particular, this implies that $\eta'(t_0) >  0$ for some $t_0 >0.$  Again by reparametrizing, we may assume that $\eta'(0) >  0.$
Moreover, since $\eta$ is non-constant and $\eta'$ is monotonically non-decreasing as $t \to \infty$, we conclude that
$$
\eta'(t) >0
$$
for $t > 0.$   Combining with \thetag {6.10}, this implies that
$$
\eta' \le \sqrt{\frac{\rho}{n-2}}\, \eta + C_0,\tag 6.11
$$
where $C_0^2 = \max\{C, 0\}.$  Dividing through by $\eta$, and integrating over the interval $(0, t)$ yield
$$
\log \eta(t) - \log \eta(0) \le \frac{r_\rho(t)}{\sqrt{n-2}} + C_0 \int_0^t \eta^{-1}
$$
with
$$
r_\rho(t) = \int_0^t \sqrt{\rho}.
$$
In particular, this implies
$$
\eta(t) \le \eta(0)\,\exp\left(C_0 \int_0^t \eta^{-1}\right)\, \exp\left(\frac{r_\rho(t)}{\sqrt{n-2}}\right).
$$
Substituting back into \thetag{6.11}, we obtain
$$
\eta'(t) \le \sqrt{\frac{\rho}{n-2}}\, \eta(0)\,\exp\left(C_0 \int_0^t \eta^{-1}\right)\, \exp\left(\frac{r_\rho(t)}{\sqrt{n-2}}\right) + C_0.\tag 6.12
$$

On the other hand, combining \thetag {6.9} and \thetag {6.11}, we have
$$
\split
\eta'' &\ge \frac{\rho}{n-2}\,\eta\\
& \ge \sqrt{\frac{\rho}{n-2}}\, (\eta' -C_0).
\endsplit\tag 6.13
$$
Let us first assume that $\eta'(t) > C_0$ as $t \to \infty.$  Since $\eta'$ is monotonically non-decreasing, there exists $t_0 >0$ such that 
$$
\eta'(t) > C_0 \ \text{for}\ t\ge t_0 \tag 6.14
$$
and we can rewrite \thetag{6.13} as
$$
(\log (\eta' - C_0))' \ge \sqrt{\frac{\rho}{n-2}}
$$
for $t \ge t_0.$
Integrating from $t_0$ to $t$, we obtain
$$
\eta'(t) - C_0 \ge (\eta'(t_0) - C_0) \exp \left( \frac{r_\rho(t) - r_\rho(t_0)}{\sqrt{n-2} }\right).
$$
Comparing this with \thetag{6.12}, we conclude that
$$
(\eta'(t_0) - C_0) \exp \left( - \frac{r_\rho(t_0)}{\sqrt{n-2}} \right) \le \sqrt{\frac{\rho}{n-2}}\, \eta(0)\,\exp\left(C_0 \int_0^t \eta^{-1}\right).
$$

Let us now assume the contrary that
$$
\liminf_{t \to \infty} \rho(t) = 0.
$$
If 
$$
\liminf_{t\to \infty} \sqrt{\rho} \,\exp\left(C_0 \int_0^t \eta^{-1}\right)=0, \tag 6.15
$$
then this will provide a contradiction.  In fact, since \thetag{6.14} implies that
$$
\split
C_0\int_0^t \eta^{-1} &< C_0 \int_0^{t_0} \eta^{-1} +C_0 \int_{t_0}^t (C_0 (t-t_0) + \eta(t_0))^{-1}\\
& \le C_1+ \log t 
\endsplit
$$
for some constant $C_1,$ \thetag {6.15} is fulfilled if
$$
\liminf_{t\to \infty} \rho(t) \,t^2 = 0.
$$

On the other hand, if 
$$
\liminf_{t\to \infty} \rho(t)\,t^2 >0,
$$
then there exists a constant $\beta >0$ such that
$$
\rho(t) \ge \beta t^{-2}\tag 6.16
$$
for sufficiently large $t>0$.  However, integrating $\eta'(t) > C_0$ for $t \ge t_0$, we have
$$
\eta(t) \ge \eta(t_0) + C_0( t - t_0).
$$
Using this and \thetag{6.16} in \thetag {6.9}, we conclude that
$$
\split
\eta'' &= (n-2)^{-1} \rho\,\eta\\
& \ge C_2 t^{-1}
\endsplit
$$
for some constant $C_2 >0$ and for sufficiently large $t.$ Hence
$$
\eta'(t) \ge C_2 \log t 
$$
 and
$$
\eta(t) \ge C_2 t\log t 
$$
for sufficiently large $t >0.$
Using this estimate on \thetag {6.9} again, we have
$$
\eta'' \ge C_3 t^{-1} \log t.
$$
Integrating again yields
$$
\eta' \ge \frac{C_3}2 (\log t)^2
$$
and
$$
\eta \ge \frac{C_3}2 t(\log t)^2
$$
for sufficiently large $t.$  Note that this estimate of $\eta$ implies that
$$
\int_0^\infty \eta^{-1} < \infty,
$$
hence the condition \thetag{6.15} is again satisfied as long as
$$
\liminf_{t\to \infty} \rho(t) = 0.
$$

Let us now consider the possibility that $\eta'(t) \le C_0$ for all $t.$  In this case, this is equivalent to taking
$$
\split
\int_0^\infty \eta'' &= \lim_{t \to \infty} \eta'(t) \\
&= C_3
\endsplit
$$
for some constant $C_3 \le C_0.$
Of course, this also implies that
$$
\eta(t) \sim C_3 t
$$
as $t \to \infty.$  Conversely, if we take $v(t)$ to be a nonnegative integrable function with
$$
\int_0^\infty v \le C_0,
$$
then one can take $\eta'' = v$ and $\eta'$ will satisfy
$$
\eta' \le C_0.
$$
For this choice, one checks rather easily that both \thetag {6.7} and \thetag {6.8} are satisfied.

With the above discussion, we conclude that 
$$
\liminf_{t\to \infty} \rho(t) > 0
$$
unless
$\eta''$ is integrable with 
$$
\split
\int_0^\infty \eta'' & = C_3\\
&\le C_0.
\endsplit
$$
In this case, $\eta$ must be asymptotically a linear function with the property
$$
\eta(t) \sim C_3 t
$$
as $t \to \infty$ for some constant $C_3 \le C_0.$ The condition that $\eta''$ is integrable is equivalent to the condition that
$$
\int_0^\infty t\,\rho(t)\,dt < \infty,
$$
which implies
$$
\liminf_{t\to \infty} t^2 \,\rho(t) = 0.
$$
Also, 
$$
\split
\int_0^\infty \sqrt{\rho} &\le \sqrt{n-2} \left(\int_0^\infty  \eta'' \right)^{\frac 12} \,\left( \int_0^\infty \eta^{-2} \right)^{\frac 12}\\
& < \infty
\endsplit
$$
 implies that the metric $ds_{\rho}^2 = \rho \, ds^2$ is incomplete. This contradicts property ($\Cal P_\rho$). \qed
\enddemo

Note that in the case when $\eta''$ is integrable, even though the $\rho$-metric is not complete, the volume growth of $M$ is polynomial of order $t^n$ 
and the finiteness theorem is still valid according to Corollary 4.2.

\heading
\S7 Parabolic Ends
\endheading

In this section, we consider the issue of number of parabolic ends.
For a 3-manifold with property ($\Cal P_\rho$),
one can also deal with the parabolic ends with the same assumption on the Ricci curvature as for the nonparabolic ends.

\proclaim{Theorem 7.1} Let $M^3$ be a complete manifold with property ($\Cal P_\rho$).  Suppse the Ricci curvature of $M$ satisfies the lower bound
$$
\text{Ric}_M (x) \ge -2 \rho(x)
$$
for all $x \in M.$  If $\rho$ satisfies the growth estimate

$$
\liminf_{R\to \infty}  (R^{-1}\,S(R))=0,
$$
then either
\roster
\item $M$ has only one end; or
\item $M$ has two nonparabolic ends and is given by $M = \Bbb R \times N$ with the warped product metric
$$
ds_M^2 = dt^2 + \eta^2(t)\,ds_N^2,
$$
for some positive function $\eta(t)$, and some compact manifold $N.$ Moreover, $\rho(t)$ is a function of $t$ alone satisfying
$$
\eta''\, \eta^{-1} = \rho.
$$
\item $M$ has one nonparabolic end and one parabolic end and is given by $M = \Bbb R \times N$ with the warped product metric
$$
ds_M^2 = dt^2 + \eta^2(t)\,ds_N^2,
$$
for some positive function $\eta(t)$, and some compact manifold $N.$ Moreover, $\rho(t)$ is a function of $t$ alone satisfying
$$
\eta''\, \eta^{-1} = \rho.
$$
\endroster
\endproclaim

\demo{Proof} According to Theorem 5.2, $M$ either has only one nonparabolic end or it must be the warped product with two nonparabolic ends.  
Hence we may assume that $M$ has one nonparabolic end and also a parabolic end $E.$ The theorem of Li-Tam \cite{L-T2}, together with a result of Nakai \cite{N}, 
asserts that one can construct a positive harmonic function $f$ with the property that
$$
\lim_{x \to \infty}f(x) = \infty \qquad \text{for} \qquad x \in E,
$$
and
$$
\liminf_{x \to \infty} f(x) = 0 \qquad \text{for} \qquad x \in M\-E.
$$

As in the case of Theorem 5.2, we consider the function $g = |\n f|^{\frac12}.$  In view of Lemma 4.1, we have
$$
\d g \ge -\rho\,g
$$
and the theorem follows by showing that
$$
\d g = - \rho\,g.
$$
Following a similar argument as in the proof of Theorem 5.2, the estimate on $M\-E$ for the term
$$
\int_{M\- E} |\n \phi|^2\, g^2
$$
can be shown to tend to 0.  We only need to deal with the term
$$
\int_E |\n \phi|^2\,g^2.
$$
Indeed, just like the case for the nonparabolic end, we just choose $\phi =\psi \, \chi$ where
$$
\chi(x) = \left\{ \aligned 1 \qquad \qquad \qquad & \qquad \text{on} \qquad \Cal L(0, T) \cap E\\
(\log T)^{-1}\,(2\log T - \log f) & \qquad \text{on} \qquad \Cal L(T, T^2)\cap E\\
0 \qquad \qquad \qquad & \qquad \text{on} \qquad \Cal L(T^2, \infty)\cap E,\endaligned \right.
$$
and
$$
\psi(x) = \left\{ \aligned 1 \qquad & \qquad \text{on} \qquad B_\rho(R-1) \cap E\\
R- r_\rho(x) & \qquad \text{on} \qquad (B_\rho(R) \- B_\rho(R-1) ) \cap E\\
0 \qquad & \qquad \text{on} \qquad  E\- B_\rho(R). 
\endaligned \right.
$$
In this case,
$$
\int_E |\n \phi|^2\, g^2 \le 2 \int_E |\n \psi|^2 \, \chi^2\, |\n f| + 2 \int_E |\n \chi|^2\, \psi^2\, |\n f|.\tag7.1
$$
Using \thetag{5.9} and a similar argument, the first term can be estimated by
$$
\split
\int_E |\n \psi|^2\, \chi^2 \, |\n f| & \le \int_{\Omega} \rho\,|\n f|\\
& \le S(R+1) \, \int_{\Omega} \rho\, f\\
& \le S(R+1)\, T^2 \int_{(B_\rho(R) \- B_\rho(R-1))\cap E} \rho,
\endsplit
$$
where 
$$
\Omega = \Cal L(0, T^2) \cap (B_\rho(R) \- B_\rho(R-1)) \cap E
$$
and
$$
S(R+1) = \sup_{B_\rho(R+1)} \sqrt{\rho}.
$$
However, combining with the estimate from Theorem 3.1 we conclude that
$$
\int_E |\n \psi|^2\, \chi^2 \, |\n f| \le C\,S(R+1)\, T^2\, \exp(-2R).\tag7.2
$$

The second term can be estimated by
$$
\split
\int_E |\n \chi|^2\, \psi^2\, |\n f| &\le
 (\log T)^{-2}\,\int_{\Cal L(T,T^2)\cap B_\rho(R) \cap E} |\n f|^3\,f^{-2}\\
&\le  (\log T)^{-2}\, S(R+1)\,\int_{\Cal L(T,T^2)\cap B_\rho(R) \cap E} |\n f|^2\,f^{-1}\\
& \le (\log T)^{-2}\, S(R+1)\,\int_{T}^{T^2} t^{-1} \int_{\ell(t)} |\n f|\, dA\, dt\\
&= C\,(\log T)^{-1}\,S(R+1) .
\endsplit
$$
Combining with \thetag{7.1} and \thetag{7.2}, we have
$$
\int_E |\n \phi|^2\, g^2 \le C\,S(R+1)\,\left( (\log T)^{-1} + T^2\, \exp(-2R) \right).
$$
As in the proof of Theorem 5.2, the right hand side tends to $0$ for a sequence of $R_i \to \infty$ by choosing $T = \exp(R\,q(R))$ with 
$q(R) =\sqrt{R^{-1}\,S(R+1)}$, and using the assumption on $S(R).$\qed
\enddemo

We now turn to the case $n\ge 4.$ 

\proclaim{Theorem 7.2} Let $M^n$ be a complete manifold of dimension $n \ge 4$ with property ($\Cal P_\rho$).  Suppose the Ricci curvature of $M$ satisfies the lower bound
$$
\text{Ric}_M (x) \ge -\frac{4}{n-1}\, \rho(x)
$$
for all $x \in M.$  If $\rho$ satisfies the property that

$$
\lim_{x\to \infty} \rho(x)=0,
$$
then $M$ has only one end.
\endproclaim

\demo{Proof} Since $\frac{4}{n-1} < \frac{n-1}{n-2}$ for $n\ge 4,$ Theorem 5.2 asserts that if $M$ has more than one nonparabolic ends then it must be 
given by the warped product and $\text{Ric}_{11} =- \frac{n-1}{n-2}\, \rho$ which is impossible.  Hence $M$ has exactly one nonparabolic end.  
Assuming that $M$ has another end $E$ that is parabolic, we construct a positive harmonic function $f$ similar to the case of Theorem 7.1.  
Again, we let $g =|\n f|^{\frac 12}.$  Lemma 4.1 asserts that it satisfies
$$
\d g \ge -\frac{2}{n-1}\,\rho\, g - \frac{n-3}{n-1}\,g^{-1}|\n g|^2.
$$
Again, because of Lemma 4.1, we will show that this inequality is indeed an equality and conclude that $M = \Bbb R \times N$ with the warped product metric
$$
ds_M^2 = dt^2 + \eta^2(t) \, ds_N^2.
$$
Moreover, the Ricci curvature in the $\frac{\p}{\p t}$ direction satisfies
$$
\split
 - \frac4{n-1}\,\rho& =
\text{Ric}_{11} \\
&= -(n-1)\eta''\, \eta^{-1}.
\endsplit \tag 7.3
$$

Following the argument of Theorem 5.2, we consider a cut-off function $\phi$ and
the integral \thetag{5.5}.  A similar argument shows that the second term on the right hand side of \thetag{5.5} can be written as
$$
\split 
2  \int_M &\phi\,g\,\la \n \phi, \n g\ra \\
&= \frac 12\int_M \la \n (\phi^2),\n (g^2)\ra\\ 
& = - \int_M  \phi^2\,g\,\d g - \int_M \phi^2\,|\n g|^2\\ 
&= \frac{2}{n-1} \int_M \phi^2\,\rho\,g^2 +\frac{n-3}{n-1} \int_M \phi^2\, |\n g|^2- \int_M \phi^2\,|\n g|^2 -\int_M \phi^2\,g\,h,
\endsplit
$$
where
$$
h = \d g+\frac{2}{n-1}\rho \,g + \frac{n-3}{n-1}\,g^{-1}|\n g|^2.
$$
Hence together with \thetag{5.5}, we have
$$
\int_M |\n (\phi\,g)|^2 + \int_M \phi^2\, g\, h = \frac{2}{n-1} \int_M \phi^2\,\rho\,g^2 + \frac{n-3}{n-1} \int_M \phi^2\, |\n g|^2 + \int_M |\n \phi|^2\, g^2.
$$
Applying \thetag{5.5} again, we obtain
$$
\int_M |\n (\phi\,g)|^2 + \frac{n-1}2 \int_M \phi^2\, g\, h =  \int_M \phi^2\,\rho\,g^2 - (n-3) \int_M \phi\,g\, \la \n \phi, \n g \ra +  \int_M |\n \phi|^2\, g^2.
$$
Together with the weighted Poincar\'e inequality, this implies that
$$
\frac{n-1}2 \int_M \phi^2\, g\, h \le -(n-3) \int_M \phi\, g\,\la \n \phi, \n g \ra + \int_M |\n \phi|^2\, g^2. \tag7.4
$$
As in the proof of Theorem 5.2, we need to choose the cut-off function $\phi$ so that the right hand side tends to 0 and the theorem follows.

Note that by a theorem of Nakai \cite{N} (also see \cite{N-R}), the positive harmonic function $f$ can be taken to be proper on the parabolic end $E$, i.e., 
$$
\liminf_{x \to \infty} f(x) = \infty
$$
for $x \in E.$
To deal with the right hand side of \thetag{7.4} on $E$,  for any $\beta >1$, we define
$$
\phi = \left\{ \aligned 1 \qquad \qquad & \qquad \text{on} \qquad \Cal L (0, T) \cap E\\
(\log \beta)^{-1} ( \log \beta T - \log f)& \qquad \text{on} \qquad \Cal L (T, \beta T) \cap E\\
0 \qquad \qquad & \qquad \text{on} \qquad \Cal L(\beta T, \infty) \cap E.\endaligned \right.
$$
Integrating by parts, we have
$$
\split
-\int_M \phi\,g\, \la \n \phi\, \n g \ra &= -\frac 14 \int_{\Cal L(T, \beta T) \cap E}  \la \n \phi^2, \n g^2 \ra\\
& = \frac 14 \int_{\Cal L(T, \beta T) \cap E} \d (\phi^2) g^2   + \frac 12 \int_{\ell (T) \cap E} \phi_\nu\, g^2,
\endsplit\tag7.5
$$
where $\nu$ is the unit normal to $\ell(t)$ given by $|\n f|\, \nu = \n f.$  Using the definition of $\phi$, we obtain
$$
\split
\int_{\Cal L(T, \beta T) \cap E} \d (\phi^2)\, g^2
&= 2\int_{\Cal L(T, \beta T) \cap E} |\n \phi|^2\, g^2 + 2 \int_{\Cal L(T, \beta T)\cap E} \phi\, \d \phi\, g^2\\
& = 2(\log \beta)^{-2} \,\int_{\Cal L(T, \beta T) \cap E} |\n f|^3\, f^{-2} \\
& \qquad + 2(\log \beta)^{-2}\,\int_{\Cal L(T, \beta T)\cap E} (\log \beta T - \log f) |\n f|^3\, f^{-2}.
\endsplit \tag 7.6
$$
Using the assumption on $\rho$ and the gradient estimate \thetag{5.9} for $f$, we conclude that there exists a constant $C >0$ such  that
$$
|\n f|(x) \le  C\,\bar S(T, \beta T)\, f(x)
$$
for all $x\in \Cal L(T, \beta T) \cap E,$ where
$$
\bar S(T, \beta T) = \sup_{\Cal B} \sqrt{\rho}
$$
with the supremum taken over the set 
$$
\Cal B = \{ x \, |\, r_\rho(x, \Cal L(T, \beta T) \cap E) \le 1\}.
$$
Hence using the fact that
$\int_{\ell(t)} |\n f|$ is a constant independent of $t,$ together with \thetag{7.6}, we obtain
$$
\split
\int_{\Cal L(T, \beta T) \cap E} &\d (\phi^2) \,g^2\\
 &\le C\, \bar S(T, \beta T)\,(\log \beta)^{-2} \,\left( \int_{\Cal L(T, \beta T)\cap E} |\n f|^2\, f^{-1} \right.\\
&\qquad + \left.\int_{\Cal L(T, \beta T)\cap E} (\log \beta T- \log f) |\n f|^2\, f^{-1}\right)\\
&= C\, \bar S(T, \beta T)\,(\log \beta)^{-2}\, \left(\int_T^{\beta T} t^{-1} + \int_T^{\beta T} (\log \beta T- \log t) t^{-1}\right)\\
& = C\,\bar S(T, \beta T)\,( (\log \beta)^{-1} + 1).
\endsplit\tag7.7
$$
On the other hand, the last term on the right hand side of \thetag{7.5} is given by
$$
\split
 \frac 12 \int_{\ell (T) \cap E}  \phi_\nu\, g^2 &\le -\frac 12 (\log \beta)^{-1}\int_{\ell(T) \cap E} f_\nu\, f^{-1}\, |\n f|\\
 & \le 0.
\endsplit
$$
Combining with \thetag{7.5} and \thetag{7.7}, we obtain
$$
-\int_M \phi\, g \, \la \n \phi, \n g \ra \le C\, \bar S(T, \beta T)\, ((\log \beta)^{-1} + 1).
$$
The second term on the right hand side of \thetag{7.4} can be estimated similarly and we have
$$
\int_M |\n \phi|^2\, g^2 \le C\, \, \bar S(T, \beta T)\, (\log \beta)^{-1}.
$$
Setting $\beta = T$ and using the assumption that $\sqrt{\rho} \to 0$, we conclude that the right hand side on \thetag{7.4} tends to 0 as $T \to \infty.$  

Note that the properness of $f$ on $E$ implies that the sublevel set $\Cal L(0 , \epsilon) \cap E = \emptyset$ for sufficiently small $\epsilon.$
By taking $M \- E$ as a nonparabolic end, we choose $\phi= \psi\, \chi$ as in Theorem 5.2.  In particular,
we set
$$
\chi(x) = \left\{ \aligned 0 \qquad \qquad \quad & \qquad \text{on} \qquad \Cal L(0, \delta \epsilon) \\
(-\log \delta)^{-1} (\log f - \log (\delta \epsilon))  & \qquad \text{on} \qquad \Cal L(\delta\epsilon, \epsilon) \\
1 \qquad \qquad \quad &  \qquad \text{on} \qquad \Cal L(\epsilon, \infty) \cap (M \- E), 
\endaligned \right.
$$
and
$$
\psi(x) = \left\{ \aligned 1\quad & \qquad \text{on} \qquad B_\rho(R-1)\\
R-r_\rho  & \qquad \text{on} \qquad B_\rho(R) \- B_\rho(R-1)\\
0 \quad & \qquad \text{on} \qquad M\-B_\rho(R) \endaligned \right.
$$
The second term on the right hand side of \thetag{7.4} on $M\- E$ can be estimated exactly as in the case when $n=3$ in the proof of Theorem 5.2.  To deal with the first term on the right hand side of \thetag{7.4}, we use
$$
-2\int_M \phi\,g\, \la \n \phi\, \n g \ra = - \int_M \psi\, \chi^2 \la \n \psi, \n (g^2) \ra - \int_M \psi^2\, \chi \la \n \chi, \n (g^2)\ra. \tag7.8
$$
The first term on the right hand side of \thetag{7.8} can be estimated by
$$
\split
- \int_M \psi\, \chi^2 \la \n \psi, \n (g^2) \ra &\le \int_M \psi \, \chi^2 |\n \psi|\, |\n (g^2)|\\
& \le \int_{\Omega} \sqrt{\rho} \,|\n (g^2)|\\
& \le \left(\int_{\Omega} \rho \right)^{\frac 12} \left( \int_{\Omega} |\n (g^2)|^2 \right)^{\frac12},
\endsplit \tag7.9
$$
where
$$
\Omega = (B_\rho(R) \- B_\rho(R-1)) \cap \Cal L(\delta \epsilon, \infty) \cap (M\- E).
$$ 
Note that the estimate \thetag{2.9} asserts that
$$
(\delta\, \epsilon)^2\, \int_{\Omega} \rho \le C\, \exp(-2R).\tag 7.10
$$
Also, since the Bochner formula (Lemma 4.1) implies that
$$
\d (g^2) \ge - \frac{4}{n-1} \rho\, g^2,
$$
if $\tau$ is a nonnegative compactly supported function, then
$$
\split
- \frac{4}{n-1} \int_M \tau^2\, \rho\, g^4 &\le \int_M \tau^2\, (g^2) \d (g^2)\\
& = -2 \int_M \tau\,g^2\, \la \n \tau, \n (g^2) \ra - \int_M \tau^2\, |\n (g^2)|^2\\
& \le 2 \int_M |\n \tau|^2\, g^4 - \frac 12 \int_M \tau^2\, |\n (g^2)|^2,
\endsplit
$$
hence
$$
\int_M \tau^2\, |\n (g^2)|^2 \le \frac{8}{n-1} \int_M \tau^2\, \rho\, g^4 + 4\int_M |\n \tau|^2\,g^4.\tag7.11
$$
Let us set
$$
\tau = \left\{ \aligned 0 \qquad & \qquad \text{on} \qquad B_\rho(R-2) \cup (M \- B_\rho(R+1))\\
r_\rho - R +2 & \qquad \text{on} \qquad B_\rho(R-1) \- B_\rho (R-2)\\
1 \qquad & \qquad \text{on} \qquad B_\rho(R) \- B_\rho(R-1)\\
R-r_\rho+1 \quad & \qquad \text{on} \qquad B_\rho(R+1)\- B_\rho(R). \endaligned \right.
$$
Then \thetag{7.11} implies
$$
\split
\int_{B_\rho(R) \- B_\rho(R-1)} |\n (g^2)|^2 &\le C\, \int_{B_\rho(R+1)\- B_\rho(R-2)} \rho\, g^4\\
&\le C\,S^2(R+1) \int_{B_\rho(R+1)\- B_\rho(R-2)} |\n f|^2, 
\endsplit \tag7.12
$$
where $S(R+1) = \sup_{B_\rho(R+1)} \sqrt{\rho}.$  Applying Corollary 2.3 to \thetag{7.12} and combining with \thetag{7.9} and \thetag{7.10}, we conclude that
$$
- \int_M \psi\, \chi^2 \la \n \psi, \n (g^2) \ra \le C\, (\delta\, \epsilon)^{-1}\,S(R+1)\, \exp(-2R).\tag7.13
$$

To estimate the second term on the right hand side of \thetag{7.8}, we integrate by parts and get
$$
\split
-\int_M \psi^2\, \chi\, \la \n \chi, \n (g^2) \ra & =- \int_{\Cal L(\delta \epsilon, \epsilon)} \psi^2\, \chi\, \la \n \chi, \n (g^2) \ra\\
&= \int_{\Cal L(\delta\epsilon, \epsilon)} \psi^2\, \chi\, \d \chi \,g^2+ \int_{\Cal L(\delta \epsilon, \epsilon) } \psi^2\, g^2\, |\n \chi|^2 \\
& \qquad + 2 \int_{\Cal L(\delta \epsilon, \epsilon)} \psi\, \chi \la \n \psi, \n \chi \ra\, g^2 - \int_{\ell(\epsilon)} \psi^2\, \chi\,\chi_\nu\, g^2 \\
& \qquad + \int_{\ell(\delta \epsilon)} \psi^2\, \chi\, \chi_\nu\, g^2,
\endsplit\tag7.14
$$
where $|\n f|\,\nu = \n f.$
Using the definition of $\chi$, the two boundary terms become
$$
\split
- \int_{\ell(\epsilon)} \psi^2\, \chi\,\chi_\nu\, g^2 + \int_{\ell(\delta \epsilon) } \psi^2\, \chi\, \chi_\nu\, g^2 
&= -(-\log \delta)^{-1}\,\int_{\ell(\epsilon)} \psi^2 f_\nu\,f^{-1}\,g^2\\
& \le 0.
\endsplit
$$
Hence \thetag{7.14} becomes
$$
\split
-\int_M \psi^2\, \chi\, \la \n \chi, \n (g^2) \ra  &\le \int_{\Cal L(\delta\epsilon, \epsilon)} \psi^2\, \chi\, \d \chi \,g^2
+2 \int_{\Cal L(\delta \epsilon, \epsilon)} \psi^2\, g^2\, |\n \chi|^2 \\
& \qquad +  \int_{\Cal L(\delta \epsilon, \epsilon)} \chi^2\, |\n \psi|^2\, g^2.
\endsplit\tag7.15
$$
We can write the term
$$
\split
\int_{\Cal L(\delta\epsilon, \epsilon)} \psi^2\, \chi\,  \d \chi \, g^2
& = -(-\log \delta)^{-2} \int_{\Cal L(\delta \epsilon, \epsilon) \cap B_\rho(R)} \psi^2\, g^2\,(\log f - \log \delta \epsilon) |\n f|^2\, f^{-2}\\
& \le 0.
\endsplit \tag 7.16
$$
As in the case of Theorem 5.2, the other term in \thetag{7.4}, and the terms in \thetag{7.15} together have the estimate given by \thetag{5.16}.  
Hence combining with \thetag{7.4}, \thetag{7.13}, \thetag{7.15} and \thetag{7.16}, we conclude that
$$
\frac{n-1}2 \int_M \phi^2\, g\, h \le C\, S(R+1)\,((\delta \epsilon)^{-1}\, \exp(-2R) + (-\log \delta)^{-1}).
$$
Since $\rho$ is bounded, by first letting $R \to \infty$ and then by setting $\delta = \epsilon \to 0,$ the right hand side must tend to 0.  
This proves that $h$ must be identically 0 and by Lemma 4.1, we conclude that $M$ is given by the warped product.  At this point, we would like to point out that the argument at the nonparabolic end only requires that
$$
\liminf_{R\to \infty} S(R+1) \, \exp(-2R) = 0.
$$

Note that using \thetag{7.3} and
$$
\text{Ric}_{\alpha \alpha} \ge -\frac 4{n-1}\, \rho,
$$
the same argument as in Theorem 6.3 asserts that
$$
\liminf_{x \to \infty} \rho(x) >0
$$
on the nonparabolic end.  However, this contradicts the assumption on $\rho$ and the theorem follows. \qed
\enddemo

\heading
\S8 Nonexistence Results for Parabolic Ends
\endheading

In this section, we will discuss some nonexistence and uniqueness results for the case when
$$
\text{Ric}_M(x) \ge - \frac 4{n-1}\, \rho(x)
$$
for a manifold with property ($\Cal P_\rho$).  The first theorem gives some indication of why we cannot prove Theorem 7.2 without the stringent assumption on $\rho.$  
The second theorem gives a restriction on the warped product situation.

\proclaim{Theorem 8.1} Let $M^n$ be a complete manifold with property ($\Cal P_\rho$) of dimension $n \ge 4.$  Suppose that the Ricci curvature of $M$ is bounded by
$$
\text{Ric}_M(x) \ge -\frac 4{n-1}\, \rho(x)
$$
for all $x \in M.$  Assume that $\rho(x)=\rho(r(x))$ depends only on the distance function $r(x)$
to a smooth compact subset $\Omega \subset M$ and it satisfies the conditions
$$
(\rho^{-\frac14})''(r) \ge 0 \quad \text{for all} \quad r \ge r_0.
$$
Then $\rho$ must be bounded from above at infinity.
\endproclaim

\demo{Proof}  In terms of Fermi coordinates $(\theta, r)$ emanating from $\p \Omega$, let us write the volume form of $M$ as
$$
dV = J(\theta, r)\, d\theta\, dr,
$$
where $d \theta$ is the volume form of $\p \Omega.$  A standard variational argument (see \cite{L1}) asserts that if we set 
$$
\psi(\theta, r) = J^{\frac 1{n-1}}(\theta, r),
$$
then $\psi$ must satisfy the differential inequality
$$
\psi'' \le \frac{4}{(n-1)^2}\,\rho\, \psi,
$$
where the derivatives are taken with respect to the $r$ variable.  Let us consider the function
$$
g(r) = C\, \rho^{-\frac 14}(r)\, \exp \left(\frac 2{n-1}r_\rho(r)\right).
$$
Direct differentiation of $g$ and  the fact that $r'_{\rho} = \sqrt{\rho}$ imply
$$
\split
g' &= -\frac{C}{4} \rho^{-\frac 54}\,\rho'\, \exp \left(\frac 2{n-1}r_\rho(r)\right) + 
\frac {2C}{n-1}\, \sqrt{\rho}\,\rho^{-\frac 14}\,  \exp \left( \frac 2{n-1}r_\rho(r) \right) \\
&= \left(- \frac 14\,\rho'\,\rho^{-1} + \frac 2{n-1}\,\sqrt{\rho} \right)\,g.
\endsplit
$$
and
$$
g'' = \left(- \frac 14 \rho''\,\rho^{-1} +\frac 5{16} (\rho')^2\, \rho^{-2} + \frac 4{(n-1)^2} \, \rho \right)g.
$$
The assumption on $\rho$ implies that 
$$
g''(r) \ge \frac 4{(n-1)^2} \, \rho \, g(r) \qquad \text{for all} \qquad r \ge r_0.
$$

Observe that  since $\rho^{\frac 12}$ is positive, there exists $r_1 \ge r_0$ such that $(\rho^{-\frac 12})'(r_1) > -\frac 4{n-1}.$  
Renaming $r_1$ to be $r_0$ if necessary, we conclude that $g(r_0) >0$ and $g'(r_0)>0.$  Hence by taking $C>0$ sufficiently large, we may assume that
$$
g(r_0) \ge \psi(\theta, r_0)
$$
and
$$
g'(r_0) > \psi'(\theta, r_0).
$$
Moreover, 
$$
(g(r) - \psi(\theta, r))''\ge \frac 4{(n-1)^2}\,\rho\, (g(r) - \psi(\theta, r))\tag 8.1
$$
for all $r \ge r_0.$  This implies that 
$$
g(r)\ge \psi(r)
$$
for all $r \ge r_0$ as otherwise $g - \psi$ will have a local positive maximum which will violate \thetag {8.1}.  
In particular,
$$
J(\theta, r) \le C\,\rho^{-\frac{n-1} 4}\, \exp (2r_\rho(r)).\tag8.2
$$

Since $M$ is nonparabolic, it must have at least one nonparabolic end $E.$  Theorem 3.1 implies that
$$
\int_{E_\rho(R+1) \- E_\rho(R)} \rho\, dV \ge C_1\, \exp(2R).
$$
Hence substituting the upper bound of the volume form \thetag{8.2} into this estimate, we obtain the inequality
$$
\int_R^{R+1} \rho^{-\frac{n-1}4 +1}\, dr \ge C_2
$$
for some constant $C_2 >0$.  Using the identity
$$
\int_R^{R+1} \rho^{\frac12} \,dr = 1,
$$
we conclude that
$$
\sup_{R\le r(r_\rho) \le R+1} \rho^{-\frac{n-3}4}(r) \ge C_2.
$$
for all $R.$  In particular, since $n \ge 4$,  there exists a sequence $r_i \to \infty$ such that
$$
\rho^{-\frac 14}(r_i) \ge C_3^{\frac 1{n-3}}.\tag8.3
$$
We now claim that
$$
\rho^{-\frac 14}(r) \ge C_3^{\frac 1{n-3}}
$$
for sufficiently large $r.$  Indeed, if this is not the case, then we can find $r_0$ and $\bar r_0$ such that there exists $r_i \in [r_0, \bar r_0]$ with
$$
\rho^{-\frac 14}(r_0) < C_3^{\frac 1{n-3}}
$$
and
$$
\rho^{-\frac 14}(\bar r_0) < C_3^{\frac 1{n-3}}.
$$
However, the assumption that $(\rho^{-\frac 14})'' \ge 0$ implies that
$$
\rho^{-\frac 14}(r) < C_3^{\frac 1{n-3}} 
$$
on the interval $ [r_0, \bar r_0]$ by the maximum principle.  This violates inequality \thetag{8.3}, hence the proposition follows. \qed
\enddemo

\proclaim{Theorem 8.2}Let $M^n = \Bbb R \times N^{n-1}$ be the warped product manifold with metric given by
$$
ds^2_M = dt^2 + \eta^2(t)\, ds^2_N.
$$
Assume that $M$ has property ($\Cal P_\rho$), and whose Ricci curvature is bounded by
$$
\text{Ric}_M(x) \ge -\frac 4{n-1}\, \rho(x)
$$
for all $x \in M.$  If we assume that $n \ge 4$ and if we denote
$$
\inf \text{Ric}_N(x) = C_N,
$$
where the infimum is taken over all points $x\in N$ and all tangent directions of $N,$  then $C_N \ge 0$.  Moreover if $C_N = 0,$ 
then $\rho$ must be identically constant given by $\rho = a^2$ for some $a >0$. In this case, $\eta(t)$ is given by either
$$
\eta(t) = \exp (at+t_0)
$$
or
$$
\eta(t) = \exp(-at + t_0)
$$
for some fixed $t_0.$ 
\endproclaim
\demo{Proof}
According to the computation in \S6, the function
$$
f(t) = \int_0^t \eta^{-(n-1)}\, ds
$$
is a harmonic function on $M.$  Moreover, a direct computation yields that if $g= |\n f|^{\frac 12}$, then
$$
\split
\d g
&=  \left(\frac{\p^2}{\p t^2} + (n-1) \eta' \, \eta^{-1} \, \frac{\p }{\p t}\right) \,\eta^{-\frac {n-1}{2}}\\
&= -\frac{n-1}2 \eta''\, \eta^{-1}\, g - \frac{(n-1)(n-3)}4 \, (\eta')^2\, \eta^{-2}\, g\\
& = \frac 12 \text{Ric}_{11}\, g -\frac{(n-1)(n-3)}4 \, (\eta')^2\, \eta^{-2}\, g.
\endsplit\tag8.4
$$
On the other hand, using \thetag{6.3}, \thetag{6.4} can be written as
$$
\split
\text{Ric}_{\alpha \alpha} &= - \eta''\, \eta^{-1}   -(n-2) (\eta')^2\, \eta^{-2} + \eta^{-2}\,  \bar{\text{Ric}}_{\alpha \alpha}\\
& = \frac 1{n-1}\, \text{Ric}_{11} - (n-2) (\eta')^2\,\eta^{-2} + \eta^{-2}\,  \bar{\text{Ric}}_{\alpha \alpha}
\endsplit
$$
for $ 2 \le \alpha \le n.$  Substituting this into \thetag{8.4}, we have
$$
\d g = \frac{(n-1)} {4(n-2)} \, \text{Ric}_{11}\,g 
+ \frac {(n-1)(n-3)}{4(n-2)}\, \text{Ric}_{\alpha \alpha}\,g -\frac {(n-1)(n-3)}{4(n-2)} \bar{\text{Ric}}_{\alpha \alpha}\, \eta^{-2}\, g.
$$
Using the lower bound of the Ricci curvature, we can write this as
$$
\d g \ge -\rho\,g -\frac {(n-1)(n-3)}{4(n-2)} \bar{\text{Ric}}_{\alpha \alpha}\, \eta^{-2}\, g .\tag8.5
$$
Observe that when $n=3$, the term $\text{Ric}_{\alpha \alpha}$ does not appear and the lower bound of the Ricci curvature has only been applied to the $\text{Ric}_{11}$ term.

By taking the infimum over all points $\bar x \in N$ and $\alpha,$ \thetag{8.5} becomes
$$
\d g \ge -\rho\,g -\frac {(n-1)(n-3)}{4(n-2)} \,C_N\, \eta^{-2}\, g .\tag8.6
$$

Let us assume that $C_N \le 0.$ Then \thetag{8.6} implies that
$$
\d g \ge -\rho\,g. \tag8.7
$$
Since $g = \eta^{-\frac {n-1}{2}}$, the integral
$$
\split
\int_{-t}^t g^2\, dV &= \int_{-t}^t \eta^{-(n-1)}\, \eta^{n-1}\,dt\\
&= 2t\\
& = o(t^2).
\endsplit
$$
Hence, we can apply the argument similar to the proof of Corollary 4.3  to conclude that \thetag{8.7} is indeed an equality.  
In particular, since $n \ge 4,$ we conclude that $C_N = 0$, 
$$
\split
-\frac 4{n-1}\, \rho&=\text{Ric}_{11} \\
&= -(n-1) \eta''\, \eta^{-1},
\endsplit
$$
and
$$
\split
-\frac 4{n-1}\, \rho&=\text{Ric}_{\alpha \alpha} \\
&= - \eta''\, \eta^{-1} - (n-2)\, (\eta')^2 \eta^{-2}.
\endsplit\tag8.8
$$
In particular, this implies that
$$
\eta''\, \eta^{-1} = (\eta')^2 \, \eta^{-2}
$$
and hence
$\log \eta$ is a linear function.  Plugging this back into \thetag {8.8}, we conclude that $\rho$ is a constant. This proves the theorem. \qed
\enddemo

\enddocument